\documentclass[reqno]{amsart}

\usepackage{amsmath}
\usepackage{amssymb}
\usepackage{graphicx}
\usepackage{color}
\usepackage[all]{xy}
\usepackage{amsthm}
\usepackage{enumerate}
\usepackage{tikz}
\usetikzlibrary{shapes.geometric}
\usepackage{overpic}
\usepackage{bbding}

\usepackage[colorlinks,citecolor=blue,linkcolor=blue]{hyperref}
\usepackage[nameinlink]{cleveref}

\theoremstyle{plain}
\newtheorem{thm}{Theorem}[section]
\newtheorem{lem}[thm]{Lemma}

\newtheorem{prop}[thm]{Proposition}
\newtheorem{cor}[thm]{Corollary}

\newtheorem{conj}[thm]{Conjecture}
\newtheorem{claim}[thm]{Claim}

\theoremstyle{definition}
\newtheorem{define}[thm]{Definition}
\newtheorem{example}[thm]{Example}
\newtheorem{rem}[thm]{Remark}
\newtheorem{conv}[thm]{Convention}

\numberwithin{equation}{section}

\makeatletter
\def\widebreve{\mathpalette\wide@breve}
\def\wide@breve#1#2{\sbox\z@{$#1#2$}%
     \mathop{\vbox{\m@th\ialign{##\crcr
\kern0.08em\brevefill#1{0.8\wd\z@}\crcr\noalign{\nointerlineskip}%
                    $\hss#1#2\hss$\crcr}}}\limits}
\def\brevefill#1#2{$\m@th\sbox\tw@{$#1($}%
  \hss\resizebox{#2}{\wd\tw@}{\rotatebox[origin=c]{90}{\upshape(}}\hss$}
\makeatletter

\newcommand{\bound}{\partial}
\newcommand{\bdy}{\partial}

\newcommand{\from}{\colon} 
\newcommand{\Stab}{\operatorname{Stab}}
\newcommand{\Isom}{\operatorname{Isom}}

\renewcommand{\setminus}{\smallsetminus}

\newcommand{\claimqed}  {\hfill \FourStarOpen \smallskip}

\newcommand{\C}{\mathbb{C}}
\renewcommand{\H}{\mathbb{H}}
\renewcommand{\P}{\mathbb{P}}
\newcommand{\Q}{\mathbb{Q}}
\newcommand{\R}{\mathbb{R}}
\newcommand{\Z}{\mathbb{Z}}

\newcommand{\TT}{\mathcal{T}}
\newcommand{\PP}{\mathcal{P}}

\newcommand{\Fp}{\mathbb{F}_p}
\newcommand{\pp}{{\mathfrak{p}}}
\newcommand{\ok}{{\mathcal{O}_k}}
\newcommand{\okp}{{\mathcal{O}_{k_{\pp}}}}
\newcommand{\oLp}{{\mathcal{O}_{L_{\pp}}}}
\newcommand{\PSL}{\mathrm{PSL}}
\newcommand{\SL}{\mathrm{SL}}
\renewcommand{\ss}{\mathbf{s}}

\newcommand{\Hth}{\mathbb{H}^3}

\newcommand{\len}{{\operatorname{len}}}
\newcommand{\smod}[1]{{\!\!\pmod{#1}}}
\newcommand{\trace}{\operatorname{tr}}

\newcommand{\PrimeSet}{\Omega}
\newcommand{\complexConj}{\bar{\tau}}

\usepackage{accents}
\newlength{\dhatheight}

\definecolor{bettergreen}{rgb}{0,0.6,0.4}
\definecolor{mutedgreen}{rgb}{.1,.75,0.15}

\definecolor{purple}{rgb}{0.4,0,0.6}

\usepackage{caption}
\usepackage{subcaption}
\usepackage{subcaption}
\begin{document}

\title{Infinitely many virtual geometric triangulations}
\author{ David Futer}
 \address{Department of Mathematics, Temple University, Philadelphia, PA 19122}
\email[]{dfuter@temple.edu}

\author{ Emily Hamilton}
 \address{Department of Mathematics,
California Polytechnic State University,
San Luis Obispo, CA 93407}
\email[]{mhamil09@calpoly.edu}

\author{ Neil R. Hoffman}
 \address{Department of Mathematics, Oklahoma State University, Stillwater, OK 74078}
\email[]{neil.r.hoffman@okstate.edu}

\makeatletter
\@namedef{subjclassname@2020}{%
  \textup{2020} Mathematics Subject Classification}
\makeatother
\date{\today}

%

\begin{abstract}
We prove that every cusped hyperbolic $3$--manifold has a finite cover admitting infinitely many geometric ideal triangulations. 
Furthermore, every long Dehn filling of one cusp in this cover admits infinitely many geometric ideal triangulations.
This cover is constructed in several stages, using results about separability of peripheral subgroups and their double cosets, in addition to a new conjugacy separability theorem that may be of independent interest.
The infinite sequence of geometric triangulations is supported in a geometric submanifold associated to one cusp, and can be organized into an infinite trivalent tree of Pachner moves.
\end{abstract}

\maketitle

\section{Introduction}

A hyperbolic $3$--manifold $M$ is called \emph{cusped} if it is noncompact and has finite volume. Every cusped $3$--manifold $M$ admits a topological ideal triangulation: that is, a decomposition into finitely many tetrahedra whose vertices have been removed, with faces identified in pairs by affine maps. A \emph{geometric ideal triangulation} is a stronger notion, where each tetrahedron is isometric to the convex hull of $4$ non-coplanar points in $\Hth$, and where the tetrahedra are glued by isometry to give the complete hyperbolic metric on $M$. See \Cref{Def:Triangulation} for precise details. The focus of this paper is on geometric triangulations.

The presence of a geometric triangulation makes the geometry of $M$ much more accessible to both practical and theoretical study. On the practical side, geometric ideal triangulations are central to the workings of the computer program SnapPy \cite{SnapPy} that computes hyperbolic structures and rigorously verifies their geometric properties 
\cite{HIKMOT16}.
On the theoretical side, Thurston's original proof of the hyperbolic Dehn filling theorem implicitly assumed the $3$--manifold at hand admits a geometric triangulation \cite{Thurston:Notes}. Similarly, Neumann and Zagier's work on volume assumes that the complement of some closed geodesic in $M$ admits a geometric triangulation \cite{NeumannZagier}.

Despite the importance of geometric triangulations, the first part of the following conjecture has been open for multiple decades. 

\begin{conj}\label{Conj:Triangulation}
Let $M$ be a (finite volume) cusped hyperbolic $3$--manifold. Then
\begin{enumerate}[$(1) $]
\item\label{Itm:OneTriang} \textup{(Folklore)} $M$ admits at least one geometric ideal triangulation.
\item\label{Itm:InfinitelyMany} $M$ admits infinitely many geometric ideal triangulations.
\end{enumerate}
\end{conj}

In the 1980s, \Cref{Conj:Triangulation}.\eqref{Itm:OneTriang} was widely believed to follow from the work of Epstein and Penner \cite{EpsteinPenner}. More precisely, the community believed that a geometric ideal polyhedral decomposition of $M$ can always be subdivided to give a geometric ideal triangulation. It took time to realize that a naive refinement of the Epstein--Penner cell decomposition does not suffice; see the discussion of coning in \Cref{Sec:Background} for a description of some of the challenges. To our knowledge, the first record of \Cref{Conj:Triangulation}.\eqref{Itm:OneTriang} in the literature is by Petronio~\cite[Conjecture 2.3]{Petronio:IdealTriangulations}, in 2000.  See also Petronio and Porti for a useful account of the history \cite{PetronioPorti}.

By contrast, \Cref{Conj:Triangulation}.\eqref{Itm:InfinitelyMany}  is new. We propose this tantalizing strengthening of the original conjecture because searching for infinite and flexible sequences of geometric triangulations might provide a pathway to finding at least one. Indeed, our main result 
can be interpreted as a proof of concept that such  a pathway exists in the context of finite covers  and Dehn filling.

Passing to covers makes both parts of \Cref{Conj:Triangulation} more amenable. Toward Part \eqref{Itm:OneTriang} of the Conjecture, Luo, Schleimer and Tillmann showed that every cusped hyperbolic manifold $M$ has a finite cover that supports a geometric triangulation  \cite{LuoSchTill}. We recall their proof strategy in \Cref{Sec:Background}, and incorporate several of their ideas in the proof of our theorems. Our main result, in the direction of \Cref{Conj:Triangulation}.\eqref{Itm:InfinitelyMany}, is the following.

\begin{thm}\label{Thm:Main}
Let $M$ be a cusped hyperbolic $3$--manifold and $A \subset M$ a horocusp. Then there is a finite cover $\widehat M \to M$, such that $A$ lifts to a cusp $\widehat A \subset \widehat M$, with the following properties:
\begin{itemize}
\item $\widehat M$ admits infinitely many geometric ideal triangulations.
\item For every sufficiently long slope $s$ on $\bdy \widehat A$, the Dehn filling $\widehat M(s)$ admits infinitely many geometric ideal triangulations.
\end{itemize} 
\end{thm}

Part of the interest of \Cref{Thm:Main}
comes from the fact that direct constructions of geometric ideal triangulations are only known in special classes of manifolds. For example, Gu\'eritaud proved that certain well-studied  triangulations of hyperbolic once-punctured torus bundles are  geometric  \cite{GueriFut}. Futer extended Gu\'eritaud's method to  hyperbolic $2$--bridge link complements \cite[Appendix]{GueriFut}. Gu\'eritaud and Schleimer proved that if $M$ is a generic multi-cusped hyperbolic manifold, 
then long Dehn fillings of $M$ will admit geometric triangulations \cite{GS:canonical}.
Ham and Purcell found geometric ideal triangulations of highly twisted link complements, by adapting Gu\'eritaud and Schleimer's construction to some especially nice triangulations of fully augmented links \cite{HamPurcell2020}.

There have also been attacks on \Cref{Conj:Triangulation}.\eqref{Itm:OneTriang} that have attempted to subdivide a geometric polyhedral decomposition into geometric ideal tetrahedra. Wada, Yamashita, and Yoshida \cite{wada1996inequality}, building on work of Yoshida \cite{yoshida1996ideal}, described a sufficient condition on the dual 1--skeleton of a polyhedral decomposition to make such a subdivision possible. Sirotkina proved that a subdivision is always possible if each $3$--cell has at most six faces \cite{Sirotkina}. Goerner proved that a subdivision is always possible if each 3--cell is a (not necessarily regular) ideal dodecahedron \cite{goerner2017geodesic}. 
Champanerkar, Kofman, and Purcell have constructed interesting examples of link complements admitting a decomposition into regular ideal bipyramids, which can then be subdivided into geometric ideal tetrahedra \cite[Theorem 3.5]{ChampanerkarKofmanPurcell}.

To our knowledge, there is only one prior paper constructing  infinitely many geometric triangulations on the same hyperbolic manifold. Dadd and Duan showed that the figure--$8$ knot complement, which decomposes into two regular ideal tetrahedra, supports infinitely many geometric triangulations \cite{DaddDuan2016}. Their proof strategy is very delicate, in that it does not extend to the figure--$8$ sister manifold, which also decomposes into two regular ideal tetrahedra.

Given a cusped manifold $M$, the \emph{topological Pachner graph} of $M$ is the graph whose vertices are isotopy classes of (topological) ideal triangulations, with edges corresponding to $2$--$3$ moves and their inverses. (See \Cref{Def:twoThree} for the definition of a $2$--$3$ move, and \Cref{Fig:Ananas} for an illustration.)
The
\emph{geometric Pachner graph} of $M$ is the induced subgraph whose vertices are geometric ideal triangulations.
The infinitely many geometric triangulations found by Dadd and Duan \cite{DaddDuan2016} are organized in the form of an infinite ray in a single component of the geometric Pachner graph of the figure--$8$ knot complement. In a generic situation, the infinitely many geometric triangulations constructed in \Cref{Thm:Main} contain an even greater amount of structure.

\begin{thm}\label{Thm:TrivalentTreeGeneric}
Let $M$ be a cusped hyperbolic $3$--manifold containing a non-rectangular cusp. Then there exists a finite cover $\widehat M \to M$ such that the geometric Pachner graph of $\widehat M$ contains a subgraph homeomorphic to an infinite trivalent tree.
\end{thm}

The hypothesis on a non-rectangular cusp can be explained as follows. As we describe in \Cref{Sec:Background}, every non-compact end of $M$ has the form $A \cong T \times [0,\infty)$, where $T$ is a torus endowed with a Euclidean metric that is well-defined up to similarity. We say that $A$ is \emph{rectangular} if the Euclidean metric on $T$ admits a rectangular fundamental domain, and \emph{non-rectangular} otherwise. By the work of Nimershiem \cite{Nimershiem}, the Euclidean structures on cusp tori of hyperbolic $3$--manifolds form a dense subset of the moduli space of $\mathcal{M}(T^2)$. Since rectangular tori represent a codimension-one slice of $\mathcal{M}(T^2)$, one can say that a generic cusped $3$--manifold satisfies the hypotheses of \Cref{Thm:TrivalentTreeGeneric}.

The infinite trivalent tree mentioned in \Cref{Thm:TrivalentTreeGeneric} can be identified with the dual $1$--skeleton of the Farey graph. See \Cref{Def:Farey} and \Cref{Fig:Farey} for a review of the Farey graph; in brief, its vertices correspond to \emph{slopes}, or simple closed curves on a torus, and to rational numbers in $\R \P^1$. The branches of the trivalent tree of \Cref{Thm:TrivalentTreeGeneric}  limit to every point of $\R \P^1$. In particular, the infinite sequence of geometric triangulations that we will construct can be chosen to approach any rational or irrational foliation on a cusp torus of $\widehat M$.
Manifolds with rectangular cusps satisfy a slightly weaker version of \Cref{Thm:TrivalentTreeGeneric}; see \Cref{Rem:RectangularFarey} for details.

\subsection{Proof strategy}
Next, we outline the main ideas in the proofs of Theorems~\ref{Thm:Main} and~\ref{Thm:TrivalentTreeGeneric}. Both proofs use the same initial setup and general strategy. Since having a non-rectangular cusp simplifies the argument considerably, \Cref{Thm:TrivalentTreeGeneric}  will be proved first.

Let $M$ be a cusped hyperbolic $3$--manifold. We will obtain geometric triangulations by subdividing the canonical (Epstein--Penner) polyhedral decompositions of covers of $M$. \Cref{Sec:Background} reviews the Epstein--Penner construction \cite{EpsteinPenner}, emphasizing the way in which the canonical polyhedral decomposition $\PP$ depends on the choice of neighborhoods of the cusps. That section also reviews the process of subdivision via \emph{coning} and lays out a sufficient condition (involving an order on the cusps) that ensures $\PP$ can be subdivided into geometric ideal tetrahedra. See \Cref{Lem:Preorder}, which is essentially due to Luo, Schleimer, and Tillmann \cite{LuoSchTill}, for details.

In \Cref{Sec:Sequence}, we describe a particular feature of the canonical polyhedral decomposition $\PP$ that occurs in the ``generic'' scenario when a manifold $M$ has multiple cusps, one cusp $A$ is chosen to be sufficiently small, and there is a unique shortest path from $A$ to the other cusps. In this situation, Gu\'eritaud and Schleimer \cite{GS:canonical} show the canonical polyhedral decomposition $\PP$ has only one or two cells poking into this cusp $A$. These cells fit together to form a submanifold called a \emph{drilled ananas} (see \Cref{Def:Ananas}). In \Cref{Lem:Ananas}, we show that a drilled ananas admits an infinite sequence of geometric ideal triangulations. When $A$ is a non-rectangular cusp, these triangulations are arranged in a trivalent tree of $2$--$3$ moves, as described in \Cref{Thm:TrivalentTreeGeneric}.

To build covers of $M$ satisfying the above-mentioned conditions, we will need to \emph{separate} certain subgroups and subsets of $\pi_1(M)$ from group elements that cause undesired coincidences. \Cref{Sec:Separability} reviews several key definitions and results about separability that are needed for our purposes. The strongest result that is needed for the proof of \Cref{Thm:TrivalentTreeGeneric} is \Cref{Thm:DoubleCosetAbelian}, due to Hamilton, Wilton, and Zalesskii \cite{HWZ:Separability}, which provides separability of  double cosets of peripheral subgroups.

With this background in hand, we can begin to construct covers. Assuming that $M$ has a non-rectangular cusp, \Cref{Sec:NonRectangular} produces a sequence of finite covers 
$\widehat M   \to \mathring M \to M$, 
with increasingly strong  properties. In particular, $\mathring M$ contains a drilled ananas, while $\widehat M$ has a polyhedral decomposition $\widehat \PP$ that can be subdivided via coning. It will follow that $\widehat M$ admits an infinite trivalent tree of geometric ideal triangulations, establishing \Cref{Thm:TrivalentTreeGeneric}.

\subsection{New separability tools}
To prove \Cref{Thm:Main}, which handles hyperbolic manifolds with rectangular cusps and provides an additional conclusion about Dehn fillings, we need stronger separability tools than what was previously available in the literature. The following new result may be of independent interest. In the theorem statement, a \emph{peripheral} subgroup of $\Gamma = \pi_1(M)$ is a subgroup coming from the inclusion of a cusp.

\begin{thm}[Conjugacy separation of peripheral cosets]\label{Thm:ConjugacySeparability}
Let $M = \Hth / \Gamma$ be a cusped hyperbolic $3$--manifold. Let $H$ and $K$ be (maximal) peripheral subgroups of $\Gamma$ corresponding to distinct cusps of $M$. Let $g \in \Gamma$ be an element such that $K$ is disjoint from every conjugate of $gH$. Then there is a homomorphism $\varphi \from \Gamma \to G$, where $G$ is a finite group, such that $\varphi(K)$ is disjoint from every conjugate of $\varphi(gH)$.
\end{thm}

 \Cref{Thm:ConjugacySeparability} has the following topological interpretation. A maximal peripheral subgroup $H \subset \Gamma$ is the stabilizer of a horoball $\widetilde B \subset \Hth$. Given $g \in \Gamma \setminus H$, the coset $gH$ is the set of all elements of $\Gamma$ that move $\widetilde B$ to $g \widetilde B$. Connecting these two horoballs is a geodesic arc $\widetilde \beta$ that projects to an arc $\beta \subset M$. We wish to find a finite cover $\widehat M \to M$ where the cusp corresponding to $K$ lifts, and where \emph{every} preimage of $\beta$ connects distinct cusps. 
\Cref{Thm:ConjugacySeparability} provides such a cover, corresponding to the subgroup
 $\widehat \Gamma = \varphi^{-1} \circ \varphi(K)$ that contains $K$ but excludes every conjugate of $gH$. 

Several precursors of \Cref{Thm:ConjugacySeparability} appear in the recent literature on $3$--manifold groups. Given a peripheral subgroup $K$ and a single element $g \in \Gamma$ that is disjoint from every conjugate of $K$, it is straightforward to find a finite quotient that witnesses this disjointness \cite[Lemma 4.5]{HWZ:Separability}. Given non-conjugate subgroups $H$ and $K$, Chagas and Zalesskii find a finite quotient of $\Gamma$ where their images are not conjugate \cite{ChagasZalesskii}.
Given a pair of non-conjugate peripheral subgroups $H$ and $K$, Wilton and Zalesskii use an argument of Hamilton to
 construct a finite quotient $\varphi \from \Gamma \to G$, such that non-trivial elements of $\varphi(H)$ and $\varphi(K)$ always lie in distinct conjugacy classes  \cite[Lemma 4.6]{WZ:DistinguishingGeometries}. 
 
 It is worth recalling the proof of the last result. First, take a hyperbolic Dehn filling $M(s)$ corresponding to a quotient $\Gamma \to \Gamma(s)$, where $K$ stays parabolic but non-trivial elements of $H$ become loxodromic. In particular, the quotient of $H$ is represented by loxodromic matrices of trace not equal to $2$. Then, take a congruence quotient of the matrix group $\Gamma(s)$ in a matrix group over a finite ring (see \Cref{Def:Congruence}), where the traces of these loxodromic matrices can still be distinguished from $2$.
Our contribution to this narrative is that we achieve even stronger separability for non-conjugate parabolic subgroups $H$ and $K$, separating the image of $K$ from the image of every conjugate of $gH$.

The proof of \Cref{Thm:ConjugacySeparability} appears in \Cref{Sec:Separability}, and uses a similar two-step method: first construct an appropriate Dehn filling, and then analyze the congruence quotients related to the Dehn filling. As part of this analysis, we apply tools from algebraic number theory, including a theorem of Hamilton
 \cite[Corollary 2.5]{hamilton2005finite} (restated below as \Cref{Prop:OrderM}), to control the traces of an entire coset $gH$.

Using the separability \Cref{Thm:ConjugacySeparability}, we prove \Cref{Thm:Main} in \Cref{Sec:RectangularDehn}. If $\mathring M$ is a cover of $M$ containing a drilled ananas $\mathring N$, as above, we use the topological interpretation of \Cref{Thm:ConjugacySeparability} to 
construct two additional covers $\widebreve M \to \overline M \to \mathring M$ where the ananas $\mathring N$ lifts but most edges of the polyhedral decomposition $\widebreve \PP$ connect distinct cusps. In particular, $\widebreve M$ has a drilled ananas $\widebreve N$ and a polyhedral decomposition $\widebreve \PP$ that can be subdivided into ideal tetrahedra via coning, which implies infinitely many geometric triangulations. Then, we build a cover $\widehat M \to \widebreve M$ where the drilled ananas $\widebreve N$ has two distinct lifts. One of these lifts supports infinitely many geometric triangulations, while the other gets filled to obtain the Dehn filling conclusion of the theorem. The opening paragraphs of \Cref{Sec:RectangularDehn} outline this construction in much greater detail.

\subsection{Acknowledgements:} 
We thank Jessica Purcell, Saul Schleimer, and Henry Segerman for helpful discussions about triangulations. We thank Ian Agol and Matthew Stover for enlightening discussions about separating peripheral double cosets. Henry Wilton, who fielded our questions at several crucial points, deserves particular gratitude.

This project was conceived when the first author visited Oklahoma State University in November 2019, and reached a mature state when the first and third authors visited the University of Arkansas for the Redbud Topology Conference in March 2020. We thank both universities for their hospitality. From there the collaboration proceeded over Zoom, with the second author joining at this virtual stage. 
 The Redbud conference stands out in our memories as a last hurrah of in-person discussion and collaboration, for a significant time to come.

During this project, Futer was partially supported by NSF grant DMS--1907708, while Hoffman was partially supported by Simons Foundation grant \#524123.

\section{Triangulations and polyhedral decompositions}\label{Sec:Background}

This section reviews some standard definitions about hyperbolic manifolds and their polyhedral decompositions and triangulations. Then, it proves \Cref{Lem:Preorder} and \Cref{Cor:PreorderTriangulation}, which will be our main ways to obtain a geometric triangulation from a polyhedral decomposition.

For the remainder of this paper, the symbol $M$ denotes a cusped, orientable hyperbolic $3$--manifold. 
(Since Theorems~\ref{Thm:Main} and~\ref{Thm:TrivalentTreeGeneric} construct finite covers, no generality is lost in assuming that $M$ is orientable.)
We will use $\widetilde M$ to denote the universal cover of $M$, which is isometric to $\Hth$. Other decorations, such as $\widehat M$ and $\overline M$, denote finite-sheeted covers of $M$.

\begin{define}\label{Def:Lift}
Let $M$ be a cusped hyperbolic manifold, and let $\widehat f \from \widehat M \to M$ be a finite cover. Let $A \subset M$ be an embedded submanifold. We say that {$A$ lifts to $\widehat M$} if the inclusion map $\iota \from A \hookrightarrow M$ lifts to an inclusion $\widehat \iota \from A \hookrightarrow \widehat M$. In this case, the image $\widehat A = \widehat \iota(A)$ is called a \emph{lift of $A$}. The lift $\widehat A$ forms only one component of $\widehat f^{-1}(A)$, and covers $A$ with degree one.

We remark that a lift is distinct from a path-lift. If $\gamma \subset M$ is a (parametrized) closed curve based at $x$, then $\gamma$ always has a path-lift $\widehat \gamma$ starting at any preimage $\widehat x \in \widehat f^{-1}(x)$. This path-lift is only a lift if it returns to $\widehat x$.
\end{define}

\begin{define}\label{Def:Triangulation}
A \emph{geometric ideal polyhedron} $P$ is the convex hull in $\Hth$ of $n \geq 4$ non-coplanar points in $\bdy \Hth$. If $n = 4$, the polyhedron is called a \emph{geometric ideal tetrahedron}, and its isometry class is determined by the cross-ratio of its $4$ vertices. The polyhedron $P$ and its boundary $\bdy P$ inherit an orientation from the embedding $P \hookrightarrow \Hth$. 

An ideal polyhedron $P$ is called an \emph{ideal pyramid} if $P$ contains an ideal vertex $v$ (called the \emph{apex}) and a unique face $F$ not incident to $v$ (called the \emph{base}). It follows that every edge of $P$ either belongs to $\bdy F$ (in which case it is called a \emph{base edge}) or connects $v$ to a vertex of $F$ (in which case it is called a \emph{lateral edge}). Every pyramid is either an ideal tetrahedron, or has a unique choice of apex and base.

A \emph{geometric ideal polyhedral decomposition} $\PP$ is a decomposition of $M$ into geometric ideal polyhedra, glued together by orientation-reversing isometries along their boundary faces. The cusps of $M$ are therefore in bijection with the equivalence classes of ideal vertices in $\PP$.  If all the cells are ideal tetrahedra, the decomposition $\PP$ is called a \emph{geometric ideal triangulation}, and denoted $\TT$. The preimage of $\PP$ in a cover $\widehat M \to M$ is denoted $\widehat \PP$, and similarly for other decorations. 
\end{define}


\begin{conv}
All triangulations and polyhedral decompositions described below are presumed to be geometric, unless specified otherwise. While there is a rich theory of topological ideal triangulations of $3$--manifolds, sometimes endowed with extra data, our focus in this paper is on geometry.
\end{conv}

\begin{define}\label{Def:Horocusp}
A \emph{(closed) horocusp} $A$ is the quotient of a closed horoball in $\H^3$ by a discrete group $G$ of parabolic isometries, where $G \cong \Z \times \Z$. Topologically, $A$ is homeomorphic to $T^2 \times [0,\infty)$ and $\bound A$ is isometric to a flat torus. The interior of $A$ is called an \emph{open horocusp}.

If $M = \H^3 / \Gamma$ is a finite-volume hyperbolic $3$--manifold, an \emph{open horocusp in $M$} is an embedded noncompact end that is isometric to an open horocusp. A \emph{(closed) horocusp in $M$} is the closure of an open horocusp in $M$. In particular, a horocusp  $A \subset M$ is homeomorphic to $T^2 \times [0,\infty)$ with a finite number of points of tangency on $T^2 \times \{0\}$ identified in pairs.

A \emph{horocusp collection} in $M$ is a union of closed horocusps $A_1, \ldots, A_n$ containing all the noncompact ends of $M$, such that the interiors of the $A_i$ are pairwise disjoint. 
 \end{define}
 
For a hyperbolic $3$--manifold $M = \H^3 / \Gamma$, we typically work with $\widetilde M = \Hth$ in the upper half-space model. The preimage of a horocusp collection in $M$ is a collection of (closed) horoballs in $\Hth$ with disjoint interiors, called a \emph{packing}. When we mention a horoball $\widetilde A$ in $\Hth$ in this context, we implicitly assume that $\widetilde A$ is one of the horoballs in the packing, meaning $\widetilde A$ covers one of the specified horocusps of $M$. 
We further conjugate $\Gamma$ in $\Isom(\Hth) \cong \PSL(2,\C)$ so that $\infty$ is a parabolic fixed point of $\Gamma$, which means that a horoball $\widetilde A$ about $\infty$ occurs in the packing. All other horoballs in the packing are tangent to points of $\C$. The packing horoballs with largest Euclidean diameter (equivalently, the horoballs closest to $\widetilde A$)  are called \emph{full-sized}.

 \begin{define}\label{Def:Orthogeodesic}
 Let $M$ be a cusped hyperbolic $3$--manifold with horocusp collection $A_1, \ldots, A_n$. An \emph{orthogeodesic}  is an immersed geodesic segment $\gamma$ that begins at $\bdy A_i$ and ends at $\bdy A_j$, such that $\gamma$ is orthogonal to $\bdy A_i$ and $\bdy A_j$ at the respective endpoints. The case $A_i = A_j$ is permitted. If $A_i$ is tangent to $A_j$, then a point of tangency is considered an orthogeodesic of length $0$. We note that an orthogeodesic is necessarily the shortest path in its homotopy class.
 
 In a similar fashion, an orthogeodesic in $\Hth$ is the shortest path between a pair of disjoint horoballs $\widetilde A, \widetilde A'$. This path is necessarily a geodesic segment that is orthogonal to $\bdy \widetilde A$ and $\bdy \widetilde A'$.  \end{define}

A collection of horocusps in a hyperbolic manifold $M$ determines a canonical decomposition of $M$ into polyhedra, as follows.

\begin{figure}
\includegraphics{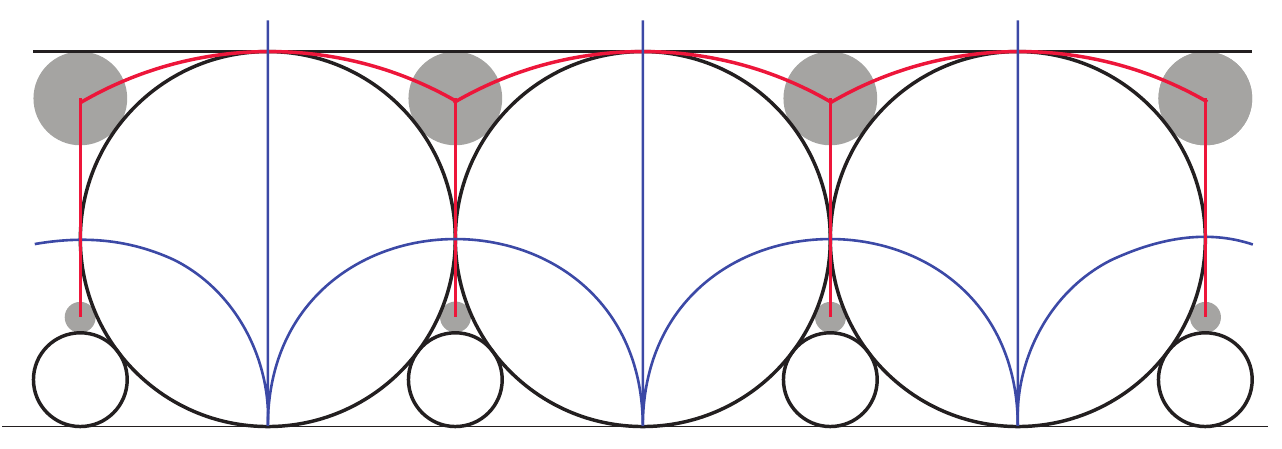}
\caption{The construction of a canonical polyhedral decomposition in a cusped hyperbolic surface. We have a horoball packing of $\H^2$ (black) and the universal cover $\widetilde \Sigma$ of the cut locus $\Sigma$ (red). Dual to $\widetilde \Sigma$ is a canonical triangulation $\widetilde \PP$ (blue). Every vertex  $\widetilde v \in \widetilde \Sigma$ is the center of a ball (grey) tangent to a maximal collection of horoballs that contain the vertices of the cell of $\widetilde \PP$ dual to $\widetilde v$.}
\label{Fig:Voronoi}
\end{figure}

\begin{define}\label{Def:CanonicalDecomp}
Let $M$ be a cusped hyperbolic $3$--manifold, endowed with a horocusp collection $A_1, \ldots, A_n$. The \emph{Ford--Voronoi domain} $\mathcal F \subset M$ consists of all points of $M$ that have a unique shortest path to the union of the $A_i$. The complement $\Sigma = M \setminus \mathcal F$, called the \emph{cut locus}, is a $2$--dimensional cell complex consisting of finitely many totally geodesic polygons. The combinatorial dual of $\Sigma$ is denoted $\mathcal P$ and called the \emph{canonical polyhedral decomposition} determined by $(M, A_1, \ldots, A_n)$. This polyhedral decomposition has one geodesic edge for each polygonal face of $\Sigma$, one totally geodesic $2$--cell for each edge of $\Sigma$, and one $3$--cell for each vertex of $\Sigma$. The edges of $\PP$ are bi-infinite extensions of orthogeodesics between the cusps. See \Cref{Fig:Voronoi} for a $2$--dimensional example.

The top-dimensional cells of $\PP$ can be characterized as follows. By construction, every $3$--cell $P \subset \PP$ is dual to a vertex $v \in \Sigma$. There is a metric ball $D$ centered at $v$, which is tangent to some number of horocusps (corresponding to the ideal vertices of $P$), and disjoint from their interiors. Furthermore, the collection of cusps tangent to $D$ is maximal with respect to inclusion.
\end{define}

In the context of closed surfaces, the construction of the canonical decomposition $\PP$ dates back to the work of Voronoi and Delaunay in the early 20th century.
Epstein and Penner  \cite{EpsteinPenner} gave a characterization of $\mathcal P$ using convexity in the hyperboloid model of $\Hth$. As a consequence, $\mathcal P$ is sometimes called the Delaunay or Epstein--Penner decomposition of $M$.

The canonical polyhedral decomposition $\PP$ determined by a choice of horocusps is always geometric. Thus every cusped hyperbolic $3$--manifold admits a geometric polyhedral decomposition. Furthermore, one may attempt to subdivide the polyhedra of $\PP$ into tetrahedra by \emph{coning}.

\begin{define}\label{Def:Coning}
Let $P$ be a (geometric) ideal polyhedron, and let $v$ be an ideal vertex of $P$. The \emph{coning of $P$ from $v$} is the decomposition of $P$ into (geometric) ideal pyramids whose apex is $v$ and whose bases are the polygonal faces of $P$ not incident to $v$. If $P$ is an ideal pyramid and $w$ is a vertex of the base of $P$, the coning of $P$ from $w$ results in ideal tetrahedra, because every face of $P$ not incident to $w$ is an ideal triangle. 
\end{define}

\begin{define}\label{Def:ConingPreorder}
Let $\PP$ be a (geometric) ideal polyhedral decomposition of $M$, and let $V$ denote the set of cusps of $M$. Consider a strict partial order $\prec$ on $V$. Observe that $\prec$ imposes a (strict) partial order on the vertices of any polyhedron $P \subset \PP$, because vertices of $P$ map to cusps of $M$. The \emph{coning of $P$ induced by $\prec$} is the following subdivision: if $P$ has a unique $\prec$--minimal vertex $v$, then $P$ is coned from $v$; otherwise, $P$ is not subdivided at all.

The \emph{iterated coning of $P$ induced by $\prec$} is the following two-step procedure. First, cone $P$ from its unique $\prec$--smallest vertex (if such a vertex exists), which either leaves $P$ unchanged or decomposes it into pyramids. Second, cone each pyramid from the unique $\prec$--smallest vertex of its base (if such a vertex exists). 
\end{define}

\begin{lem}\label{Lem:Preorder}
Let $M$ be a cusped hyperbolic $3$--manifold, and $\prec$ a strict partial order on the set of cusps of $M$. Let $\PP$ be a (geometric) ideal decomposition of $M$, with the property that every polyhedron $P \subset \PP$ has a unique $\prec$--minimal vertex $v_P$. Then the iterated coning of $\PP$ induced by $\prec$ produces a well-defined subdivision of $\PP$ into geometric ideal pyramids.

Furthermore, if $\prec$ gives a total order of the vertices of every polyhedron, then the iterated coning of $\PP$ produces a geometric ideal triangulation.
\end{lem}

\begin{proof}
Suppose $P$ and $P'$ are polyhedra of $\PP$ that are identified along a face $F$. We need to check that the iterated coning of $P$ and $P'$ induces the same subdivision of the face $F$. We show this by considering two cases.

\textbf{Case 1:}  $F$ does not have a unique $\prec$--minimal vertex. In this case, we claim that $F$ will not be subdivided at all. For, polyhedron $P$ will be coned from its unique minimal vertex $v_P$, which is not contained in $F$ by hypothesis. This produces a collection of pyramids, with $F$ a base of one of the pyramids. Since $F$ does not have a unique minimal vertex, the second stage of iterated coning does not subdivide $F$ at all. An identical argument applies to $P'$, proving the claim.

\textbf{Case 2:} $F$ has a unique $\prec$--minimal vertex $w$. In this case, we claim that $F$ will be subdivided by coning from $w$.
If $w = v_P$ is minimal in all of $P$, then $P$ will be subdivided into pyramids by coning from $w$, hence $F$ will be also. Otherwise, if $w \neq v_P$, then $P$ will be subdivided into pyramids by coning from $v_P$, and $F$ will be the base of one of these pyramids. At the second stage of the iterated coning, the pyramid in $P$ containing $F$ will be coned from the minimal vertex of $F$, namely $w$. An identical argument applies to $P'$, proving the claim.

Finally, observe that if $\prec$ gives a total order of the vertices of each cell, then we must be in Case 2: every face $F$ has a minimal vertex. Thus every face is subdivided into triangles, and every pyramid is subdivided into geometric ideal tetrahedra.
\end{proof}

The ``furthermore'' statement of \Cref{Lem:Preorder} was previously observed 
 by Luo, Schleimer, and Tillmann  \cite[Lemma 7]{LuoSchTill}. They also used the separability of peripheral subgroups (\Cref{Prop:PeripheralSeparability}) to show that every cusped $3$--manifold $M$ with a geometric polyhedral decomposition $\PP$ has a finite cover $\widehat M$ such that any order on the cusps of $\widehat M$ imposes a total order on the vertices of each polyhedron of $\widehat \PP$. Compare \Cref{Lem:NoDiagonal} below. Consequently, $\widehat M$ has a geometric ideal triangulation.

In fact, a total order is not necessary to produce a geometric triangulation:
 
\begin{cor}\label{Cor:PreorderTriangulation}
Let $M$ be a cusped hyperbolic $3$--manifold, and $\prec$ a strict partial order on the set of cusps of $M$. Let $\PP$ be a (geometric) ideal decomposition of $M$, with the property that every polyhedron $P \subset \PP$ has a unique $\prec$--minimal vertex $v_P$. Let $\PP'$ be the pyramidal refinement of $\PP$ guaranteed by \Cref{Lem:Preorder}. Then every choice of diagonals in the non-triangular faces of $\PP'$ leads to a decomposition of $\PP'$ into geometric ideal tetrahedra.
\end{cor}
 
 \begin{proof}
Following \Cref{Lem:Preorder}, let $\PP'$ be the subdivision into ideal pyramids coming from the iterated coning of $\PP$ induced by $\prec$. Then every non-triangular $2$--cell $F \subset \PP'$ must be the base of exactly two pyramids. Thus the non-tetrahedral pyramids of $\PP'$ are glued in pairs, with each pair forming a bipyramid that is joined to other cells along ideal triangles only. 
This means we have complete freedom to choose diagonals of every non-triangular $2$--cell $F$, subdividing the two pyramids adjacent to $F$ into tetrahedra, without impacting the choices anywhere else in the manifold. Every such choice produces a subdivision of $\PP'$ into geometric ideal tetrahedra.
 \end{proof}

 \section{An infinite tree of triangulations}\label{Sec:Sequence}
 
In the last section, we described a construction of Luo, Schleimer, and Tillmann for decomposing an ideal polyhedral decomposition $\PP$ into a geometric ideal triangulation. Having one geometric triangulation is clearly a prerequisite to having infinitely many. In this section, we describe a particular geometric feature called a \emph{drilled ananas} (see \Cref{Def:Ananas}) that admits an infinite sequence of geometric triangulations. By embedding a drilled ananas inside a triangulation of a cusped manifold $M$, we obtain an infinite sequence of ideal triangulations of $M$. {See \Cref{Lem:Ananas} for the construction of an infinite sequence of triangulations, and \Cref{Prop:FareyAnanas} for a more refined description of an infinite trivalent tree of triangulations.}


Consider a polyhedral decomposition $\PP$ of $M$. Let $A \subset M$ be a horocusp, chosen small enough that for every polyhedron $P \subset \PP$, the intersection $P \cap A$ consists of neighborhoods of ideal vertices. Then $\PP$ induces a decomposition of the torus $\bdy A$ into Euclidean polygons, which truncate the ideal vertices of polyhedra of $\PP$. We call this decomposition the \emph{cusp cellulation} of $\bdy A$, and denote it $\mathcal{C}(A)$. If $\PP$ is the canonical polyhedral decomposition (determined by some choice of cusps), then $\mathcal{C}(A)$ satisfies the \emph{Delaunay condition}: the vertices of every polygon can be inscribed on a circle, where the interior of the circle does not contain any other vertices.

In the following proposition, $A$ is a horocusp of $M$. Let $A^t$ be a sub-horocusp of $A$ such that $d(\bound A, \bound A^t) = t$. A particular feature occurs when $t$ becomes sufficiently large.

\begin{prop}[Gu\'eritaud--Schleimer \cite{GS:canonical}]\label{Prop:GS}
Let $M$ be a cusped hyperbolic $3$--manifold, endowed with a choice of horocusps $A, B_1, \ldots, B_n$ for $n \geq 1$. Assume that an orthogeodesic $\alpha$ from $A$ to $B_1$ is the  unique shortest path from $A$ to $\cup_{j=1}^n B_j$. Then, for every sufficiently small sub-horocusp $A^t \subset A$, the canonical decomposition $\PP$ determined by $A^t, B_1, \ldots,B _n$ contains a unique edge from $A^t$ to $\cup_{j=1}^n B_j$. This edge is the bi-infinite extension of $\alpha$.

Furthermore, there are one or two $3$--cells of $\PP$ that meet the cusp $A^t$. Each such $3$--cell is an ideal pyramid with an apex at $A^t$ and all lateral edges identified to $\alpha$. If $A$ is a rectangular cusp, then the single $3$--cell meeting $A^t$ is a rectangular pyramid and the induced cellulation of $\bdy A^t$ is a rectangle. Otherwise, if $A$ is a non-rectangular cusp, then the two $3$--cells meeting $A^t$ are isometric ideal tetrahedra, and
 the induced cellulation of $\bdy A^t$ consists of two isometric, acute triangles.
\end{prop}

This result is due to Gu\'eritaud and Schleimer \cite[Section 4.1]{GS:canonical}, and appears in the form of a discussion with the explicit hypothesis that the cusp torus $\bound A$ is not rectangular. For completeness, we reproduce an expanded version of their proof.

\begin{proof}
Set $M = \Hth/\Gamma$, where we view $\Hth$ in the upper half-space model. As described in \Cref{Sec:Background}, we conjugate $\Gamma$ in $\operatorname{Isom}(\Hth)$ so that a horoball $\widetilde A$ about $\infty$ covers the horocusp $A$. Then $K = \Stab_\Gamma(\infty) \cong \Z^2$ can be identified with $\pi_1 (A) \subset \Gamma$. Every preimage of $B_j$  is a horoball tangent to a point of $\C$.

By hypothesis, there is a unique shortest orthogeodesic from $A$ to $\cup_{j=1}^n B_j$, which leads from $A$ to $B_1$. Extend $\alpha$ to be a bi-infinite geodesic. After shrinking $A$ by an appropriate distance, we may further assume that $\alpha$ is the unique shortest orthogeodesic from $A$ to $A \cup (\cup_{j=1}^n B_j)$. Consequently, there is a choice of horoball $\widetilde B_1$ covering $B_1$, such that all of the full-sized horoballs tangent to points of $\C$ are in the $K$--orbit of $\widetilde B_1$.

\begin{figure}
\includegraphics[width=5in]{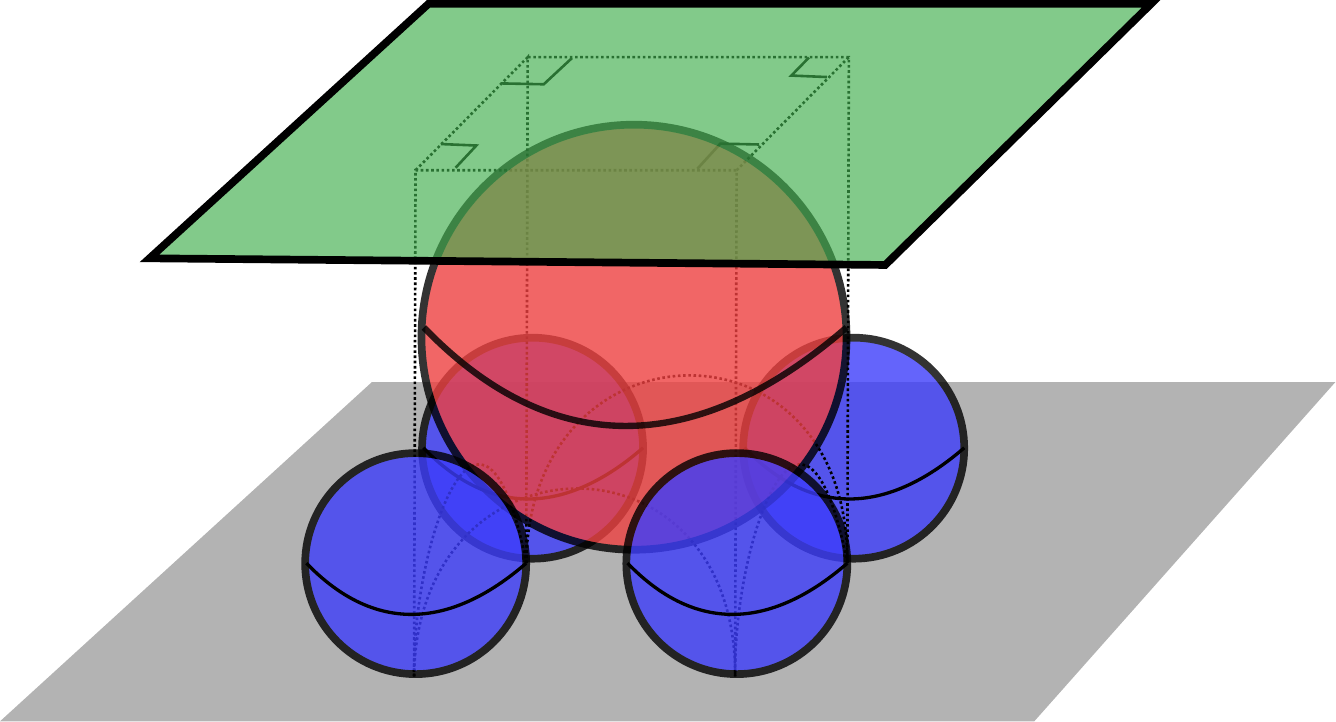}
\caption{A ball $D$ (shown in red) resting on four full-sized horoballs (blue) and tangent to a horoball about $\infty$ (green). In the setting of \Cref{Prop:GS}, all of the full-sized horoballs are in the same orbit of $K \cong \Z^2$. It follows that the ball $D$ can be tangent to four full-sized horoballs at once if and only if $K$ has a rectangular fundamental domain, as in this figure. The center of $D$ corresponds to an ideal rectangular pyramid (dotted).}
\label{fig:horoball_intro}
\end{figure}

Consider a ball $D \subset \Hth$ that rests on the collection of horoballs tangent to points of $\C$ (see \Cref{fig:horoball_intro}).
If the Euclidean diameter of $D$ is sufficiently large, then $D$ will only touch the full-sized horoballs. Furthermore, if we shrink $A$ by a sufficiently large distance $t \geq 0 $, producing a sub-horocusp $A^t$ whose preimage horoball $\widetilde A^t \subset \widetilde A$ is at sufficient Euclidean height, then $D$ will also be disjoint from $\widetilde A^t$. 
Inflating $D$ to a maximal (hyperbolic) radius produces a ball $D^+$ that is tangent to $\widetilde A^t$ and some number of horoballs from the orbit $K \cdot \widetilde B_1$, and is disjoint from all other horoballs. Observe that $D^+$ is tangent to either $3$ or $4$ horoballs in $K \cdot \widetilde B_1$, and that the case of $4$ occurs precisely when $\bdy A^t = \bdy \widetilde A^t / K$ has a rectangular fundamental domain.

Now, consider the cut locus $\Sigma^t$ corresponding to the horocusp collection $A^t, B_1, \allowbreak \ldots, B_n$. Let $\widetilde \Sigma^t$ be the preimage of $\Sigma^t$ in $\Hth$, and consider the polyhedral decomposition $\PP = \PP^t$ dual to $\Sigma^t$. By \Cref{Def:CanonicalDecomp}, every vertex $v \in \widetilde \Sigma^t$ is equidistant from some collection of horoballs (which is maximal with respect to inclusion), and conversely every point equidistant from a maximal collection of horoballs is a vertex of $\widetilde \Sigma^t$. In particular, the hyperbolic center of the ball $D^+$ constructed in the previous paragraph must be a vertex $w \in \widetilde \Sigma^t$. By \Cref{Def:CanonicalDecomp}, this vertex $w$ is dual to a polyhedron $P_w$ in the canonical decomposition $\PP$, whose ideal vertices lie in the horoballs tangent to $D^+$.
Recall that $D^+$ is tangent to $\widetilde A^t$ and either $3$ or $4$ full-sized horoballs  in $K \cdot \widetilde B_1$. 

If there are $3$ full-sized horoballs tangent to $D^+$, then $P_w$ is an ideal tetrahedron. Furthermore, up to the action of $K = \Stab_\Gamma(\infty)$, there must be exactly one other ideal tetrahedron in $\PP$ that intersects $\widetilde A^t$. In this case, the induced cusp cellulation of $\bound A^t$ consists of two isometric triangles. 
{These triangles must be acute: otherwise, the circle that circumscribes the $3$ vertices of an obtuse triangle would contain the fourth vertex of the parallelogram in its interior, violating the Delaunay condition.}
Since the two triangles of $\bdy A^t$ are isometric, the two tetrahedra meeting $A^t$ are also isometric. Furthermore, each of the two tetrahedra has $3$ edges identified to $\alpha$, corresponding to the fact that the cusp cellulation has a single vertex at $\alpha \cap \bdy A^t$.

If there are $4$ full-sized horoballs tangent to $D^+$, as in \Cref{fig:horoball_intro}, then the polyhedron $P_w$ is an ideal rectangular-based pyramid. In this case, $P_w$ is the only cell of $\PP$ meeting $A^t$, and the induced cusp cellulation of $A^t$ is a single rectangle with a vertex at $\alpha \cap \bdy A^t$. Consequently, all lateral edges of $P_w$ are identified to $\alpha$.
\end{proof}

\begin{rem}\label{Rem:Akiyoshi}
Akiyoshi  \cite{akiyoshi_finiteness} proved that the canonical polyhedral decomposition $\PP^t$ determined by the horocusp collection $A^t, B_1, \allowbreak \ldots, B_n$ must stabilize  as $t \to \infty$. Thus, for sufficiently large $t$, there is a stable geometric decomposition $\PP = \PP^t$ independent of $t$. This is our motivation for dropping the superscript $t$.
\end{rem}

In the polyhedral decomposition $\PP$ that occurs in \Cref{Prop:GS}, the cells that enter the special cusp $A^t$ fit together to form a geometric object that we call a \emph{drilled ananas}.

\begin{define}\label{Def:Ananas}
A \emph{drilled ananas} is a $3$--manifold $N$ homeomorphic to $T^2 \times [0, \infty) \setminus \{x\}$, where $x \in T^2 \times \{ 0\}$, and endowed with a complete hyperbolic metric 
with the following properties. 
For some $y > 0$, the non-compact end $T^2 \times [y, \infty) \subset N$ is isometric to a horocusp, such that each cross-section  $T^2 \times \{ y' \}$ for $y'>y$ is a horotorus. The boundary $\bound N = T^2 \times \{ 0\} \setminus \{x\}$ is made up of two totally geodesic ideal triangles, with vertices at $x$. These two triangles are glued by isometry along their edges to form a standard two-triangle triangulation of a once-punctured torus, with shearing and bending allowed along the edges of this boundary triangulation. If the base of the drilled ananas is comprised of a single totally geodesic ideal rectangle, then we make a choice of diagonal to decompose $N$ into two ideal tetrahedra. 

The horocusp $T^2 \times [y, \infty) \subset N$ is called the \emph{cusp} of $N$, while a regular neighborhood of $x \in T^2 \times \{ 0\}$ is called the \emph{thorn} of $N$.
\end{define}

Here are a few notes on terminology and past usage. The term \emph{thorn} was coined by Baker and Cooper \cite{Baker-Cooper:QFS}. A drilled ananas is a special case of a \emph{topological ideal polyhedron} in the work of Gu\'eritaud \cite{Gueritaud:SolidTori}, and a slightly less special case of an \emph{ideal torihedron} in the work of Champanerkar, Kofman, and Purcell  \cite[Definition 2.1]{ChampanerkarKofmanPurcell}. Both of these generalizations capture the idea of placing a hyperbolic structure on a $3$--manifold endowed with a polyhedral graph on its boundary. Gu\'eritaud coined the term \emph{ananas} (French for \emph{pineapple}) to describe a topological polyhedron with the topology of a solid torus. 
(Compare the definition of a filled ananas immediately above \Cref{Claim:FilledAnanas}.) The object $N$ in \Cref{Def:Ananas} can be obtained by removing the core of a solid torus, hence \emph{drilled ananas}.

With the above definition, the conclusion of \Cref{Prop:GS} can be rephrased as follows.

\begin{cor}\label{Cor:EmbeddedAnanas}
Let $M$ be a cusped hyperbolic manifold satisfying the hypotheses of \Cref{Prop:GS}, and let $\PP = \PP^t$ be the polyhedral decomposition produced by that proposition. Let $N \subset M$ be the submanifold obtained by gluing together all cells of $\PP$ that have an ideal vertex in cusp $A$, along their shared faces. Then $N$ is a drilled ananas comprised of two acute ideal tetrahedra or one ideal rectangular pyramid. 
Furthermore, $N$ is convex, with an angle less than $\pi$ at every edge of $\PP \cap \bdy N$.
\end{cor}

\begin{figure}[t!]
\begin{subfigure}[t]{.4\textwidth}
  \centering
  \includegraphics[width=.8\linewidth]{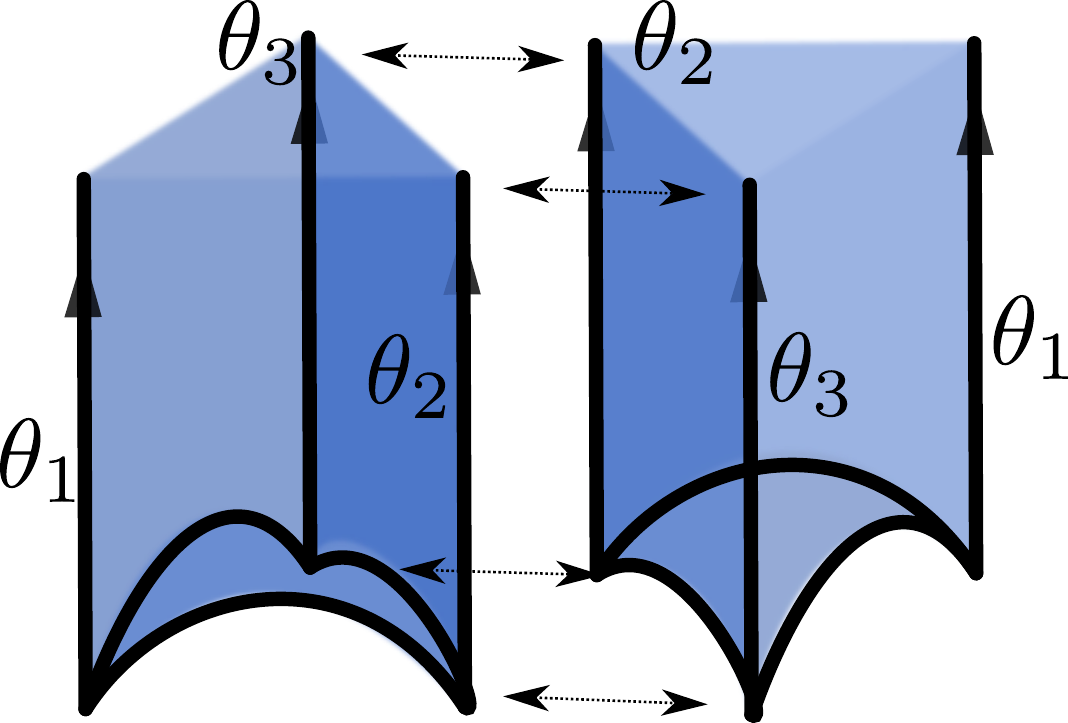}  
  \caption{Gluing two ideal tetrahedra to form a drilled ananas $N_i$.}
  \label{fig:sub-ananas-two-tets}
\end{subfigure}
\hspace{10mm}
\begin{subfigure}[t]{.4\textwidth}
  \centering
  \includegraphics[width=1.1\linewidth]{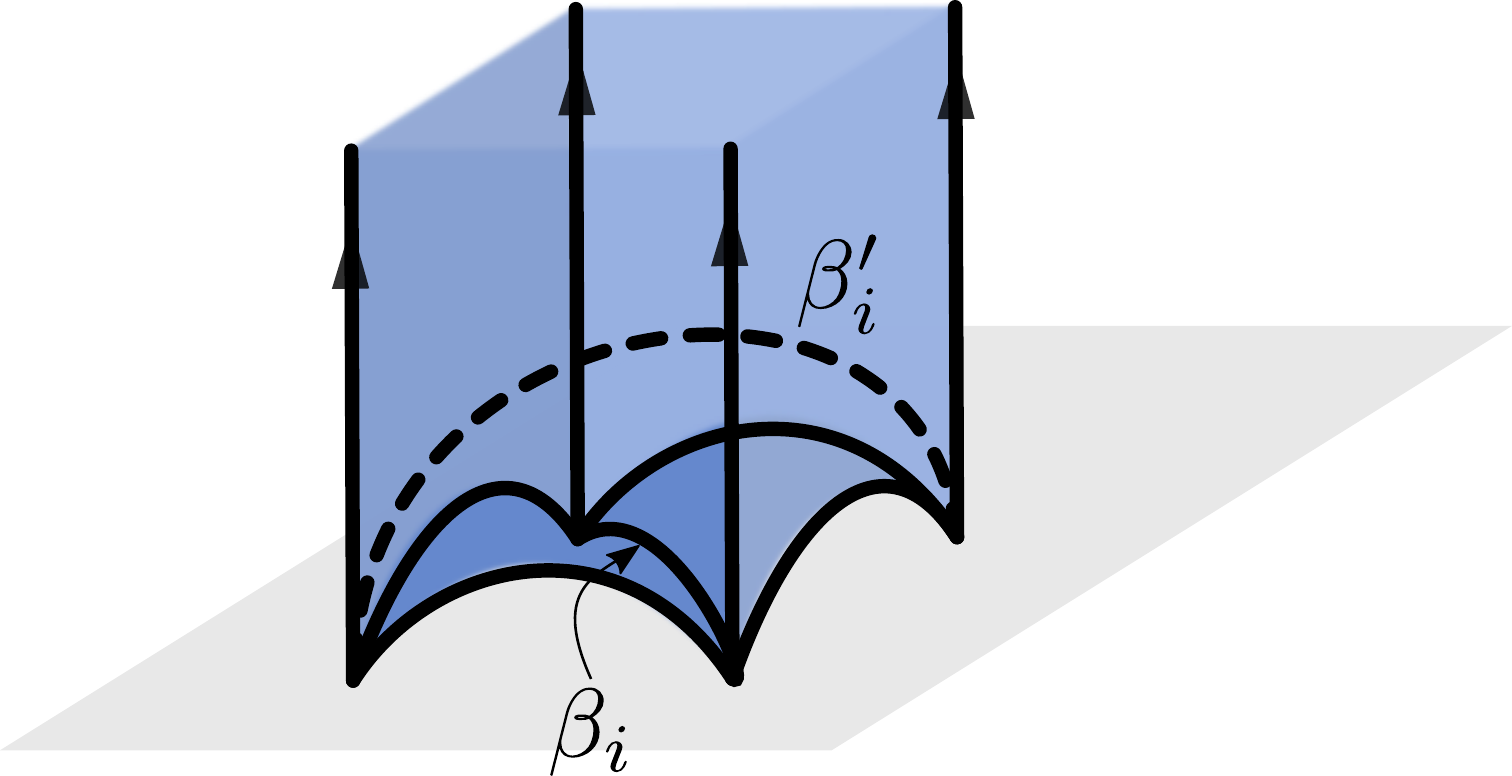}  
  \caption{In a $2$--$3$--move, the two tetrahedra glued along a vertical face, shown in part I, are replaced by three tetrahedra that share an edge dual to that face (dashed). The result appears in part III. }
  \label{fig:sub-ananas-two-tets-glued}
\end{subfigure}
\begin{subfigure}{.4\textwidth}
  \centering
\hspace{-20mm}
  \includegraphics[width=1.1\linewidth]{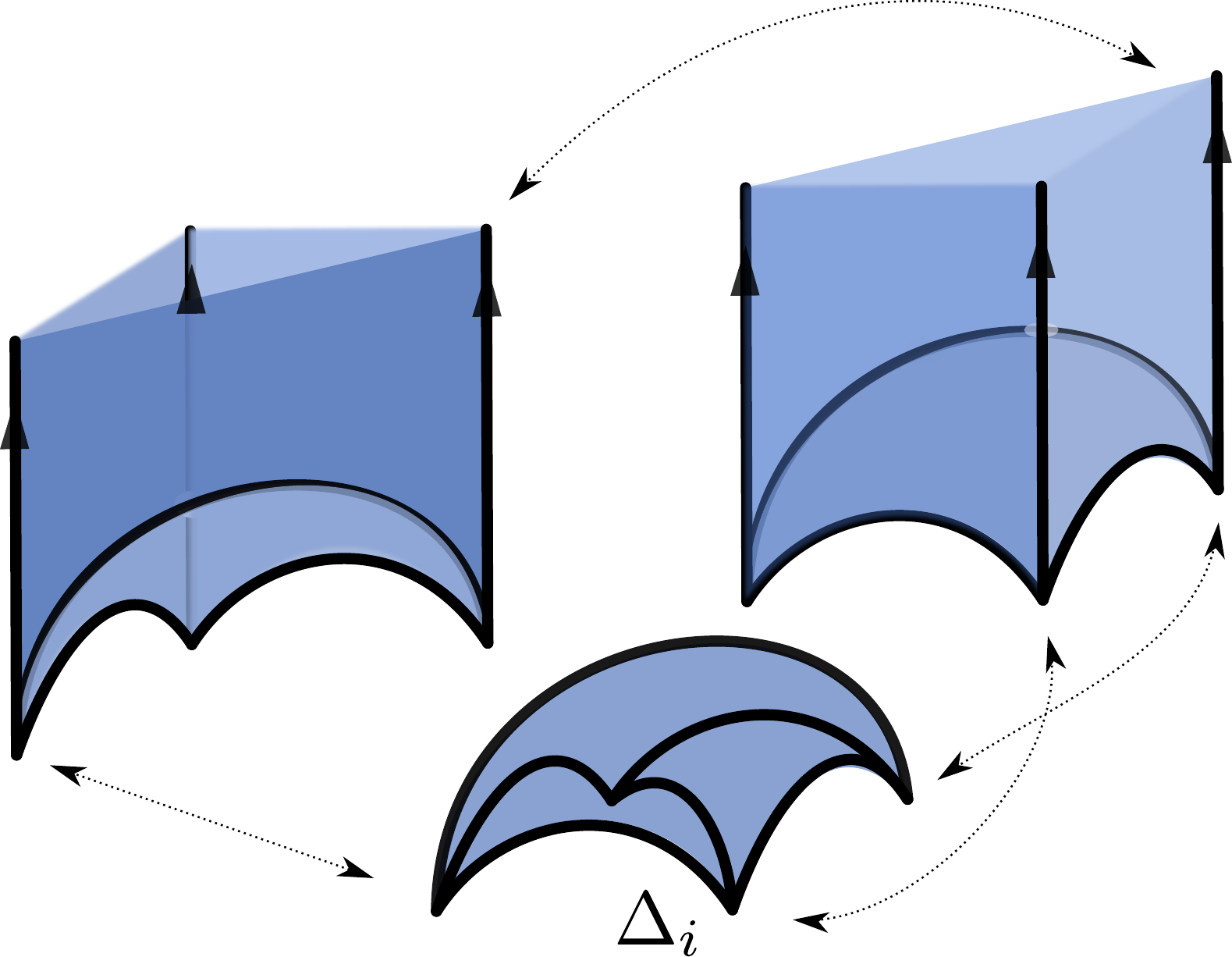}  
  \caption{After a $2$--$3$ move, the original ananas $N_i$ decomposes into three ideal tetrahedra.}
  \label{fig:sub-ananas-three-tets}
\end{subfigure} 
\hspace{10 mm}
\begin{subfigure}{.4\textwidth}
  \centering
  \includegraphics[width=1.1\linewidth]{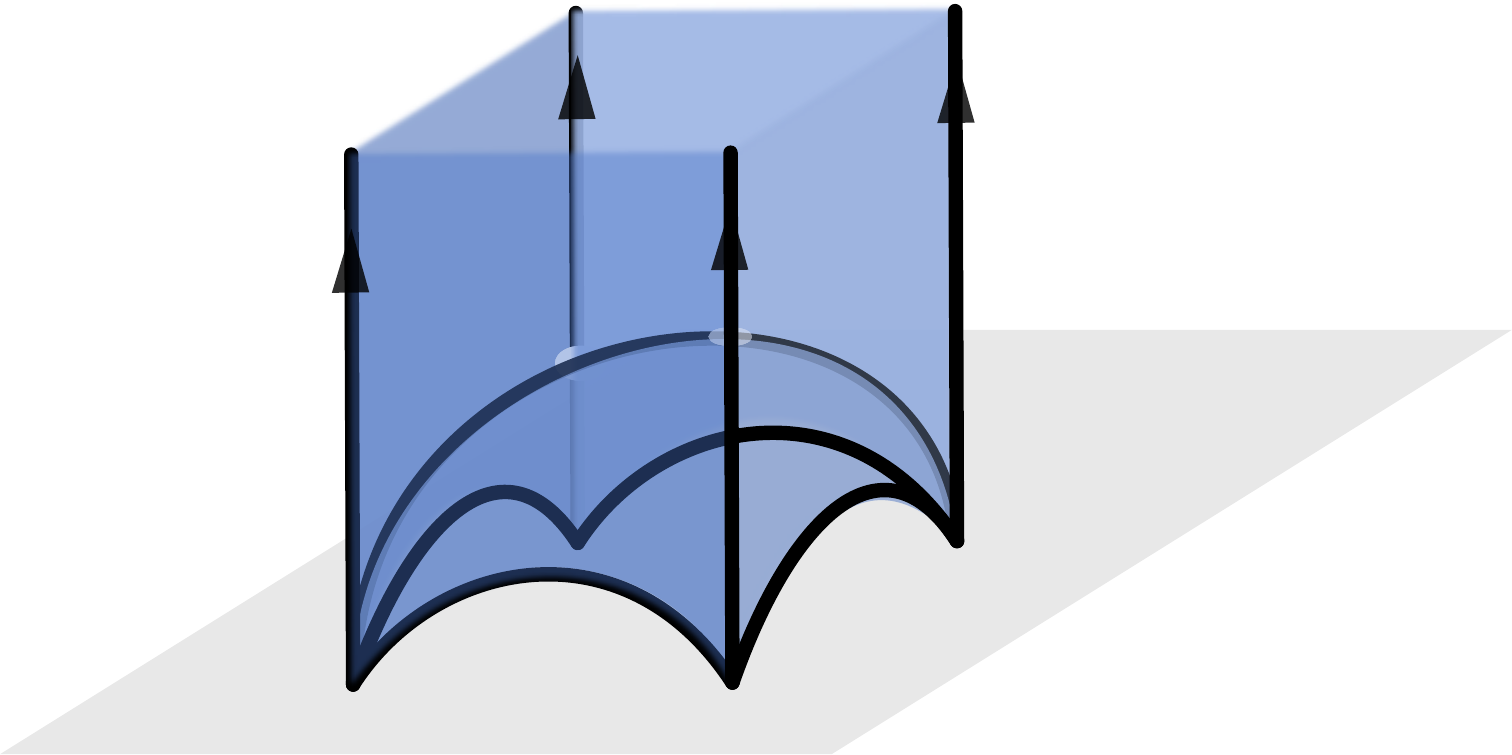}  
  \caption{The new ananas $N_{i+1}$ is the union of the two tetrahedra in $N_i$ which are incident to $\infty$. The next $2$--$3$ move will be along a face that is exterior to the fundamental domain shown here.}
  \label{fig:sub-newly-formed}
\end{subfigure}
\begin{subfigure}{.4\textwidth}
  \centering
   \end{subfigure}

\caption{The inductive step of \Cref{Lem:Ananas}. }
\label{Fig:Ananas}
\end{figure}

\begin{proof}
The conclusion that $N$ is a drilled ananas is immediate from \Cref{Prop:GS} and \Cref{Def:Ananas}. 
It remains to check that $N$ is convex.

By \Cref{Prop:GS}, each cell of $\PP$ comprising $N$ is an ideal pyramid with a base along $\bdy N$. If $N$ contains a single rectangular pyramid, we subdivide it into two isometric ideal tetrahedra $T, T'$ by choosing a diagonal along $\bdy N$. Otherwise, $N$ already consists of two isometric, acute-angled ideal tetrahedra $T, T'$.
The three lateral faces of $T$ are glued to the three lateral faces of $T'$, and all lateral edges of $T,T'$ are identified to the single geodesic $\alpha$ that connects the thorn of $N$ to the cusp of $N$.

Let $\theta_1, \theta_2, \theta_3$ be the dihedral angles of $T$ at the edges identified to $\alpha$. Since $T$ and $T'$ are necessarily isometric, these are also the dihedral angles of $T'$ at the edges identified to $\alpha$. See \Cref{Fig:Ananas}(I). Then the three internal angles along the edges of $\bdy N$ are $2\theta_1, 2\theta_2, 2 \theta_3$. 

If $T$ and $T'$ are cells of $\PP$, then \Cref{Prop:GS} says that all of their angles are acute. Thus $2 \theta_i < \pi$ for every $i$, hence $\bdy N$ is locally convex at every edge on its boundary. Otherwise, if $T$ and $T'$ were created by subdividing an ideal pyramid, the cusp $A \subset N$ is rectangular, hence the Euclidean triangles truncating the tips of $T$ and $T'$ have a right angle $\theta_1$ and two acute angles $\theta_2, \theta_3$. Consequently, the internal angles along $\bdy N$ are $2 \theta_1 = \pi$ and $2 \theta_2, 2 \theta_3 < \pi$. In either case, $\bdy N$ is locally convex, hence $N$ is convex.
\end{proof}

A key feature of a drilled ananas is that it is made up of two tetrahedra glued along three faces. Each of these three faces supports a local move, called a $2$--$3$ move.

\begin{define}\label{Def:twoThree}
Let $\TT_0$ be a topological ideal triangulation of a $3$--manifold $N$, possibly with boundary. Let $T$ and $T'$ be distinct tetrahedra in $\TT_0$ that are glued together along a face $F$, and observe that $T \cup_F T'$ is a bipyramid. A \emph{topological $2$--$3$ move} replaces $T \cup T'$  with three tetrahedra glued together along a central edge dual to $F$, while all other tetrahedra of $\TT_0$ remain the same. The resulting triangulation is denoted $\TT_1$. We say that a $2$--$3$ move is \emph{geometric} if $\TT_0$ and $\TT_1$ are both geometric triangulations. Assuming $\TT$ is geometric, a $2$--$3$ move will be geometric whenever the bipyramid $T \cup_F T'$ is strictly convex.

In \Cref{Fig:Ananas}, tetrahedra $T$ and $T'$ are shown in panel (I), with face $F$ a darker shade of blue. The dual edge to $F$ is dashed in panel (II), and the resulting three tetrahedra are shown in panel (III).
\end{define}

\begin{lem}\label{Lem:Ananas}
A drilled ananas $N$ admits an infinite sequence of geometric triangulations, connected by geometric $2$--$3$ moves. 
\end{lem}

The inductive construction that proves the lemma is illustrated in \Cref{Fig:Ananas}.

\begin{proof}
Let $N$ be a drilled ananas. By definition, $\bdy N$ is subdivided into two totally geodesic ideal triangles. Consequently, $N$ itself can be subdivided into two geometric ideal tetrahedra $T, T'$ by coning those triangles to the non-compact end at $\infty$. As in \Cref{Prop:GS}, we think of each of $T, T'$ as a triangular pyramid with a base along $\bdy N$. The three lateral faces of $T$ are glued to the three lateral faces of $T'$, and all lateral edges of $T,T'$ are identified to a single geodesic $\alpha$.

As in \Cref{Cor:EmbeddedAnanas}, let $\theta_1, \theta_2, \theta_3$ be the dihedral angles of $T$ and $T'$ at the edges identified to $\alpha$. See \Cref{Fig:Ananas}(I). Since these angles are positive and  $\theta_1 + \theta_2 + \theta_3 = \pi$, two of the three angles must be strictly less than $\pi/2$. Then the three internal angles along the edges of $\bdy N$ are $2\theta_1, 2\theta_2, 2 \theta_3$.  At least two of these angles (say, $2\theta_1$ and 
$2\theta_2$) are strictly less than $\pi$.

We will prove the lemma by induction. Set $N_0 = N$. The key inductive claim is:

\begin{claim}\label{Claim:Induct}
Let $N_i$ be a drilled ananas with a two-tetrahedron geometric triangulation. Let $\beta_i \subset \bdy N_i$ be a boundary edge such that the internal angle at $\beta_i$ is less than $\pi$. Then performing a diagonal exchange on $\beta_i$ results in a geometric ideal tetrahedron $\Delta_i$ and a new sub-ananas $N_{i+1} \subset N_i$, where $N_i = N_{i+1} \cup \Delta_i$ and $N_{i+1}$ again has a two-tetrahedron geometric triangulation.
\end{claim}

The proof of the claim is almost immediate. Removing $\beta_i$ from a triangulation of $\bdy N_i$ results in a quadrilateral. Let $\beta'_i$ be the opposite diagonal of this quadrilateral. Observe that the geodesic representative of $\beta'_i$ lies strictly inside $N_i$, because $N_i$ is locally convex at $\beta_i$. 
 The join of $\beta_i, \beta'_i$ is a (geometric) ideal tetrahedron $\Delta_i$. Removing the interior of $\Delta_i$ from $N_i$ produces a sub-ananas $N_{i+1}$ whose boundary is pleated along $\beta'_i$ and the two remaining edges of $\bdy N_i$. \claimqed

Observe that the two-tetrahedron triangulation of $N_i$ gives rise to a three-tetrahedron triangulation of $N_i = \Delta_i \cup N_{i+1}$, in a geometric $2$--$3$ move. See \Cref{Fig:Ananas}. We can now apply this move to $N_{i+1}$, and so on, resulting in an infinite sequence of geometric triangulations of $N = N_0$.
\end{proof}

\begin{cor}\label{Cor:AnanasCover}
Suppose $N$ is a drilled ananas, and $f \from \widehat N \to N$ is a finite cover. Then $\widehat N$ admits an infinite sequence of $f$--equivariant geometric triangulations.
\end{cor}

The corollary follows immediately from \Cref{Lem:Ananas} because every triangulation of $N$ lifts to a triangulation of $\widehat N$. 
We also remark that any finite cover $\widehat N \to N$ is regular and has abelian deck group, because $\pi_1(N) \cong \Z^2$ is abelian.

\subsection{Connections to the Farey complex}
We can now describe some additional structure in the set of triangulations of a drilled ananas $N$.

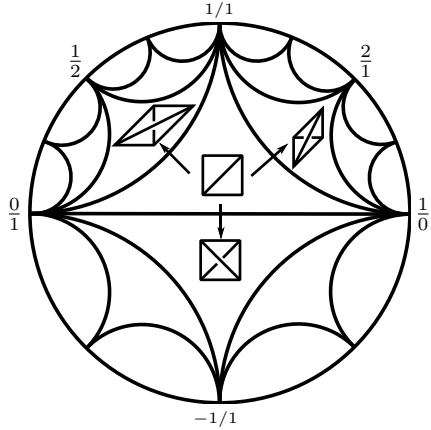
\begin{figure}
\begin{overpic}[width=2in]{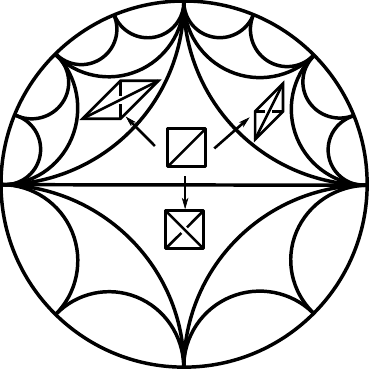}
\put(101,48){$\tfrac{1}{0}$}
\put(-6,48){$\tfrac{0}{1}$}
\put(46,102){\tiny{$1/1$}}
\put(43,-5){\tiny{$-1/1$}}
\put(86,88){$\tfrac{2}{1}$}
\put(10,88){$\tfrac{1}{2}$}
\end{overpic}
\caption{The Farey complex $\mathcal{F}$. Edges of the dual tree correspond to diagonal exchanges in a torus with one marked point. Every non-backtracking path in the dual tree, starting from the central triangle, can be realized via geometric $2$--$3$ moves in a drilled ananas. Figure from Ham and Purcell \cite[Figure 3.1]{HamPurcell2020}. }
\label{Fig:Farey}
\end{figure}

\begin{define}\label{Def:Farey}
The \emph{Farey complex} $\mathcal{F}$ is a simplicial complex whose vertices are isotopy classes of arcs on a torus $T^2$ based at a marked point $x$, and whose edges correspond to arcs that are disjoint (except at $x$). The vertices of $\mathcal F$ are commonly identified with the rational points  $\Q\P^1 \subset \R\P^1$, as follows. Endow $T^2$ with a standard Euclidean metric, with fundamental domain a unit square. Then every loop in $T^2$ based at $x$ can be pulled tight to a Euclidean geodesic of  some well-defined  slope $\Q \cup \{ \infty \}$. Conversely, every rational slope defines a unique isotopy class of arc from $x$ to $x$.

Triangles in $\mathcal{F}$ correspond to (isotopy classes of) one-vertex triangulations of $T^2$ with the vertex at $x$, or equivalently to ideal triangulations of $T^2 \setminus \{x\}$. The dual $1$--skeleton of $\mathcal{F}$ is a trivalent tree, with every edge of the dual tree corresponding to a diagonal exchange. See \Cref{Fig:Farey}.
\end{define}

Using \Cref{Def:Farey}, we can state the following stronger formulation of \Cref{Lem:Ananas}.

\begin{prop}\label{Prop:FareyAnanas}
Let $N$ be a drilled ananas with a geometric triangulation consisting of two acute-angled tetrahedra. Then $N$ admits an infinite trivalent tree of geometric $2$--$3$ moves, with vertices of the tree in natural bijection with triangles of the Farey complex $\mathcal{F}$.
\end{prop}

\begin{proof}
The proof amounts to adding some book-keeping to the proof of \Cref{Lem:Ananas}. Let $N_0 = N$, subdivided into ideal tetrahedra $T, T'$. By hypothesis, all dihedral angles of $T$ and $T'$ are acute. 
This triangulation of $N$ defines an induced cellulation of the boundary of the horocusp $ A \subset N$. In fact, this cusp cellulation is a one-vertex triangulation: the one vertex is the intersection between $\bdy A$ and the single edge $\alpha \subset N$ into cusp $A$, while the two triangles are the cross-sections $\bdy A \cap T$ and $\bdy A \cap T'$. 
Let $\delta_0$ be this triangulation of $\bdy A$. Fix a framing of $\bdy A$ so that the three slopes occurring in $\tau_0$ are $0/1$, $1/0$, and $1/1$, as in the central triangle of \Cref{Fig:Farey}.

Observe that any one-vertex triangulation of $\bdy A$ defines an ideal triangulation of $\bdy N$, by projecting outward from the cusp. In the opposite direction, any ideal triangulation of the boundary $\bdy N_i$ for some sub-ananas $N_i \subset N$ defines a one-vertex cusp triangulation. We will pass freely between the two viewpoints.

Now, let $\beta_0$ be a boundary edge of $N_0$. In order to carry out the construction of \Cref{Claim:Induct}, the internal angle of $\beta_0$ needs to be less than $\pi$. But since all dihedral angles of $T$ and $T'$ are acute, we may choose $\beta_0$ at will. Now, the new sub-ananas $N_1 \subset N_0$, constructed as in  \Cref{Claim:Induct}, will induce a new cusp cellulation $\tau_1$ of $\bdy A$, which differs from $\tau_0$ via a diagonal exchange. Thus we may choose $\tau_1$ to be any one of the three triangulations adjacent to $\tau_0$ in \Cref{Fig:Farey}.

Continuing inductively, suppose that we have constructed the sub-ananas $N_i \subset N_{i-1} \subset \ldots \subset N_0$, and that $\bdy N_i$ has ideal triangulation $\tau_i$. Then, as noted in the proof of \Cref{Lem:Ananas},
 there are two edges of $\bdy N_i$ that have interior angle less than $\pi$. The third edge necessarily has angle more than $\pi$; it is the edge $\beta'_{i-1}$ that was just created in constructing $N_i$. (See \Cref{Fig:Ananas}(IV), with indices shifted by $1$.) Thus we may choose $\beta_i$ to be any edge of $\bdy N_i$ other than $\beta'_{i-1}$, which means the new cusp triangulation $\tau_{i+1}$ is allowed to be any of the two neighbors of $\tau_i$ that are distinct from $\tau_{i-1}$. In summary, the path $\tau_0, \tau_1, \ldots$ of cusp cellulations associated with $N_0, N_1, \ldots$ is allowed to be any non-backtracking path starting from  $\tau_0$.
\end{proof} 

In any path $\tau_0, \tau_1, \ldots$ of cusp cellulations constructed in the above proof, the slopes of edges approach some limiting value in $\R\P^1$, and the edges themselves approach a foliation with the limiting slope.  Thus the geometric retriangulations of $N$ can be chosen to limit to any foliation of the torus.

The strategy for proving \Cref{Thm:Main} and \Cref{Thm:TrivalentTreeGeneric} can now come into view. Given a cusped $3$--manifold $M$, we will find a cover $\mathring M$ with a polyhedral decomposition $\mathring \PP$ as in \Cref{Prop:GS}, where two ideal tetrahedra fit together to form a drilled ananas. A further cover $\widehat M$ ensures that the other cells of $\widehat \PP$ can be subdivided into ideal tetrahedra as well. Producing these covers requires tools from subgroup separability, described in the next section.

\section{Separability}\label{Sec:Separability}

This section begins with a review of some standard results about separable subsets and subgroups in a group $G$. The main content of the section is a proof of~\Cref{Thm:ConjugacySeparability}.

\begin{define}\label{Def:Profinite}
Let $G$ be a group. The \emph{profinite topology} on $G$ is the topology whose basic open sets are cosets of finite-index normal subgroups. Since every coset of a finite-index subgroup $H \lhd G$ is the complement of finitely many other cosets of $H$, the basic open sets are also closed.

A subset $S \subset G$ is called \emph{separable} if it is closed in the profinite topology on $G$.
\end{define}


The following characterization is standard.

\begin{lem}\label{Lem:Separable}
Let $G$ be a group and $S \subset G$ a subset. The following are equivalent:
\begin{enumerate}
\item\label{Itm:Separable} $S$ is separable.
\item\label{Itm:SeparateOne} For every element $g \in G \setminus S$, there is a homomorphism $\varphi \from G \to F$, where $F$ is a finite group and $\varphi(g) \notin \varphi(S)$.
\end{enumerate}
\end{lem}

\begin{proof}
The set $S$ is closed if and only if every $g \in G \setminus S$ is contained in a basic open set disjoint from $S$.
But, by \Cref{Def:Profinite}, a basic open set containing $g$ is precisely the preimage of an element under a homomorphism to a finite group.
\end{proof}


If $M$ is a cusped hyperbolic $3$--manifold, a subgroup of $\pi_1(M)$ coming from the inclusion of a horocusp $A$ is called \emph{peripheral}. This peripheral subgroup of $\pi_1(M)$ is also \emph{maximal abelian}: it is not contained in any larger abelian subgroup.

 We will need to separate peripheral subgroups and their double cosets. 
The following separability result has been known since the 1980s, if not earlier. See e.g.\ Long \cite[page 484]{Long:TotallyGeodesic} for a proof.

\begin{prop}\label{Prop:PeripheralSeparability}
Let $M = \H^3 / \Gamma$ be a cusped hyperbolic $3$--manifold.
Then maximal abelian subgroups of $\Gamma$ are separable. In particular, peripheral subgroups of $\Gamma$ are separable.
\qed
\end{prop}

We will also need to separate peripheral double cosets, using the following theorem of Hamilton, Wilton, and Zalesskii \cite[Theorem 1.4]{HWZ:Separability}.

\begin{thm}[Hamilton--Wilton--Zalesskii \cite{HWZ:Separability}]\label{Thm:DoubleCosetAbelian}
Let $M = \H^3 / \Gamma$ be a finite-volume hyperbolic $3$--manifold. Let $H$ and $K$ be abelian subgroups of $\Gamma$. Then, for every $g \in \Gamma$, the double coset $HgK = \{ hgk : h \in H, k \in K \}$ is separable in $\Gamma$. \qed
\end{thm}

\Cref{Thm:DoubleCosetAbelian} is the strongest separability tool needed in the proof of \Cref{Thm:TrivalentTreeGeneric}.  The reader who is mainly interested in that result is invited to proceed directly to \Cref{Sec:NonRectangular}.

\subsection{Algebraic tools for separability}
The separability results that we use in this paper, including  \Cref{Thm:DoubleCosetAbelian} and \Cref{Thm:ConjugacySeparability}, are proved using tools from algebraic number theory. To set up the proof of \Cref{Thm:ConjugacySeparability}, we review some needed definitions and results.

A \emph{number field} is a finite field extension of $\Q$. An extremely useful connection between number fields and $3$--manifolds 
comes from the following result.

\begin{thm}[Thurston~\cite{Thurston:Notes}, Bass~\cite{bass1980groups}, Culler--Shalen~\cite{Culler:Lifting,CuSh}]
\label{Thm:ThurstonLiftConjugate}
Let $M = \Hth/\Gamma$ be a cusped hyperbolic $3$--manifold. Then 
\begin{enumerate}[$(1) $]
\item\label{Itm:RepLifts} $\Gamma \subset \PSL(2,\C)$ can be lifted to $\SL(2,\C)$. 
\item\label{Itm:ConjugateRing} $\Gamma$ can be conjugated in $\SL(2,\C)$ to lie in 
$\SL(2,R) \subset \SL(2,k)$, where $R$ is a finitely generated subring of a number field $k$. \qed
\end{enumerate}
\end{thm}

Conclusion~\eqref{Itm:RepLifts} was observed by Thurston~\cite[page 98]{Thurston:Notes} and carefully written down by Culler and Shalen~\cite[Proposition 3.1.1]{CuSh}, with a simplified proof by Culler~\cite[Corollary 2.2]{Culler:Lifting}.
Conclusion~\eqref{Itm:ConjugateRing} was observed by Thurston~\cite[Proposition 6.7.4]{Thurston:Notes} as an algebraic consequence of Mostow rigidity, and independently proved by Bass~\cite{bass1980groups}.

The arithmetic data associated to a ring $R$ can be used to construct finite quotients  of both rings and groups. The construction uses the following notions. 

\begin{define}\label{Def:Congruence}
If $k$ is a number field, let $\ok$ denote the ring of integers of $k$.
If $\pp$ is a non-zero prime ideal of $\ok$, we let $k_{\pp}$ denote  the $\pp$-adic completion of $k$ and 
$\okp$ 
the ring of integers of $k_{\pp}$.   The ring $\okp$ has a unique maximal ideal.  
The quotient of $\okp$ by this maximal ideal is a finite field called the \emph{residue class field} of $\okp$.  The quotient map
is called the \emph{residue class field map} with respect to $\pp$. If $R$ is a finitely generated ring in a number field $k$, then $R \subset \okp$ for all but finitely
many primes $\pp$ of $\ok$.  Restricting the residue class field map to $R$ yields a finite quotient of $R$.

Now, suppose $\Gamma \subset \SL(2, R)$, where $R$ is a ring. Let $I \subset R$ be an ideal, such that $S = R/I$ is finite.
Then the homomorphism
\[
\SL(2, R) \to \SL(2, S) 
\]
where coefficients are reduced modulo $I$ is called a \emph{congruence quotient} of PSL$(2,R)$. We also call the composition
\[
\varphi \from \Gamma \hookrightarrow \SL(2, R) \to \SL(2, S) 
\]
a \emph{congruence quotient} of $\Gamma$. (This is slightly abusive, because $\Gamma$ may fail to surject $\SL(2, S)$.)
\end{define}

We  
now describe
 three algebraic results that will be used in the proof of 
\Cref{Thm:ConjugacySeparability}. The first of these is  \cite[Theorem 2.6]{HWZ:Separability}.

\begin{prop}[Hamilton--Wilton--Zalesskii \cite{HWZ:Separability}]
\label{Prop:AddRingSep}
Let $R$ be a finitely generated ring in a number field $k$.
By fixing a $\Q$ embedding of $k$ into $\C$, we may view $k \subset \C$.
Let $\omega \in R$, and set
$Z_\omega = \{ m + n\omega \mid m, n \in \mathbb{Z} \}$ and $Q_\omega = \{ m + n\omega \mid m, n \in \Q \}$.
If $y \in R - Q_{\omega}$, then there exist a finite ring $S$
and a ring homomorphism $\rho \from R \rightarrow S$ such that
$\rho(y) \notin \rho(Z_{\omega})$. \qed
\end{prop}

The following technical lemma is a variant of \Cref{Prop:AddRingSep}.

\begin{lem}\label{Lem:LinearIndep}
Let $R$ be a finitely generated ring in a number field $k$.  By fixing a $\Q$ embedding of $k$ into $\C$,
we may view $k \subset \C$.  Let $\omega \in R - \R$. Then there is an infinite collection  $\PrimeSet$ of primes of $\Q$, such that for each prime $p \in \PrimeSet$, there exist a finite field $F_{\pp}$ 
of characteristic $p$ and a ring homomorphism $\eta_p\from R \rightarrow F_{\pp}$, such that $\{ 1, \eta_p(\omega) \}$ is linearly independent over $\Fp$. 

Consequently, $\eta_p$ has the following property. Consider an element $y_*  \in Q_\omega =\{ m + n\omega \mid m, n \in \Q \}$. Express $y_*$ in lowest terms:
\[
y_* = \frac{m_* + n_* \omega}{v_*} \quad 
\]
where $m_*, n_* \in \mathbb{Z}$ and $v_* \in \mathbb{N}$. 
If $\eta_p(y_*) = \eta_p(m + n \omega)$ for some $m, n \in \mathbb{Z}$, then $v_* m \equiv m_* \smod p$ and $v_* n \equiv n_* \smod p$.
\end{lem}

\begin{proof}
Since $\omega \in \C - \R$, the set $\{ 1, \omega \}$ is linearly independent over $\Q$.  
By a standard argument used in the proof of \Cref{Thm:DoubleCosetAbelian} (compare~\cite[page 278]{HWZ:Separability}), we can preserve this property in a finite quotient.  We include the details for completeness.
Let $L$ denote the normal closure of $k$ over $\Q$ and let $\complexConj
\in \mathrm{Gal}(L/\Q)$ represent complex conjugation.
Since $\omega \in \C - \R$, we know $\omega$ is not fixed by $\complexConj$. 
By the Tchebotarev Density Theorem,
there exist infinitely many primes $p$ of $\Q$ with unramified
extension $\pp$ in $L$ such that $\complexConj$
is the Frobenius automorphism for ${\pp} / p$. After eliminating a finite set of primes, if necessary, we may assume that
$R \subset \oLp$, where $\oLp$ denotes the ring of integers in the $\pp$-adic field $L_{\pp}$.
Given such a prime $p$, let $F_{\pp}$ denote the residue class field of $\oLp$ and let
$\Fp$ denote the finite field of $p$ elements.
Let $\eta_p$ be the composition of the inclusion map of
$R$ into $\oLp$ with the residue map:
$$ \eta_p\from R \hookrightarrow \oLp \rightarrow F_{\pp}.$$
Since $\complexConj$ is the Frobenius automorphism of $L/\Q$
with respect to $\pp / p$, $\mathrm{Gal}(L_{\pp} / \Q_{p})
= \langle \complexConj^{\prime} \rangle$ where $\complexConj^{\prime} = \complexConj$ on $L$.
Since $\complexConj(\omega) \neq \omega$, it follows that $\omega \notin \Q_{p}$.
The Galois group of $F_{\pp} / \Fp$ is also induced by $\complexConj$.
This implies that $\eta_p(\omega) \notin \Fp$ and therefore, the set $\{ 1, \eta_p(\omega) \}$ is linearly independent over $\Fp$. 

Next, suppose that $\eta_p(y_*) = \eta_p(m + n \omega)$ for some $m, n \in \mathbb{Z}$. Then 
\[\eta_p(m_* + n_* \omega) = \eta_p(v_* y_*) 
= \eta_p(v_* m + v_* n \omega).
\]  
The above equality can be rewritten as $\eta_p(v_* m - m_*) + \eta_p(v_* n - n_*) \eta_p(\omega) = 0$.  Since 
the set $\{ 1, \eta_p(\omega) \}$ is linearly independent over $\Fp$, it follows that 
$v_* m \equiv m_* \smod p$ and $v_* n \equiv n_* \smod p$.
\end{proof}

The third preliminary algebraic result is a combination of 
 \cite[Theorem 2.3 and Corollary 2.5]{hamilton2005finite}.

\begin{prop}[Hamilton \cite{hamilton2005finite}]\label{Prop:OrderM}
Let $R$ be a finitely generated ring
in a number field $k$, let $\lambda$ be
a non-zero element of $R$
that is not a root of unity, and let
$x_1, x_2, \ldots , x_j$ be non-zero
elements of $R$.
Then, for every sufficiently large integer $q$,
there exist a non-zero prime ideal $\pp$  of $\ok$, a finite field $F_\pp$, 
and a ring homomorphism $\sigma\from R \rightarrow F_\pp$
such that $R \subset \okp$, the multiplicative order of $\sigma(\lambda)$ is
equal to $q$, and $\sigma(x_i) \neq 0$, for each
$1 \leq i \leq j$. The field $F_\pp$ is the residue class field of $\okp$ and the map $\sigma$ is the restriction to $R$ of the residue class field map with respect to $\pp$. \qed
\end{prop}

\subsection{Conjugacy separation of peripheral cosets}
We can now prove \Cref{Thm:ConjugacySeparability}.

\begin{proof}[Proof of \Cref{Thm:ConjugacySeparability}]
 Let $H$ and $K$ be maximal parabolic subgroups of $\Gamma$ corresponding to 
distinct cusps of $M$, and let $g \in \Gamma$ be an element such that $K$ is disjoint from every conjugate of $gH$.
In particular, this implies $g \notin H$.
Fix a non-trivial element $h_0 \in H$.  
Since $H$ is a maximal abelian subgroup of $\Gamma$, and $g \notin H$, the commutator
$[g,h_0] = gh_0g^{-1}h_0^{-1}$ is nontrivial.  Let $A$ be the cusp of $M$ corresponding to $K$, let  $B$ be the cusp corresponding to $H$, and let $C_1, \dots, C_\ell$ be the remaining cusps. We will leave the cusp $A$ unfilled, and 
will fill the remaining cusps. So long as a tuple of slopes $\ss$ on $B, C_1, \dots, C_\ell$ avoids finitely many slopes on each cusp, 
the Dehn filled manifold $M(\ss)$ will be hyperbolic. Thus, by Thurston's hyperbolic Dehn surgery theorem, we can choose generators $h_1$ and $h_2$ of $H$ that can be completed to tuples $\ss_1$ and $\ss_2$ where $M(\ss_1)$ and $M(\ss_2)$ are both hyperbolic. Then, for $j=1,2$, the fundamental group $\pi_1(M(\ss_j))$ has a discrete, faithful representation to a group of isometries $\Gamma(\ss_j) \subset \PSL(2,\C)$. By \Cref{Thm:ThurstonLiftConjugate}, we view $\Gamma(\ss_j)$ as a subgroup of $\SL(2,\C)$.
Let $\psi_j \from \Gamma \to \Gamma(\ss_j)$ be the quotient homomorphism induced by the inclusion $M \hookrightarrow M(\ss_j)$. 
By choosing sufficiently long Dehn fillings, we may assume that 
$\psi_1([g, h_0])$ and $\psi_2([g, h_0])$ are non-trivial.  These choices ensure the following properties for $j \in \{ 1, 2 \}$:
\begin{itemize}
\item $\psi_j(K)$ is a parabolic subgroup of $\Gamma(\ss_j)$, 
\item $\psi_j(H)$
is a loxodromic subgroup of $\Gamma(\ss_j)$, 
\item $\psi_j(g) \notin \psi_j(H)$.
\end{itemize}

Before working with the Dehn filled manifolds, we examine the coset $gH$ in $\Gamma$.
By \Cref{Thm:ThurstonLiftConjugate},
we can conjugate $\Gamma$ to lie in
$\SL(2,k)$ for
some number field $k$. 
After possibly expanding $k$, and then conjugating $\Gamma$ in $\SL(2, k)$, 
we may assume that 
\[
 g = \begin{pmatrix} a & b \\ c & d \\ \end{pmatrix}  , \quad 
h_1 = \pm \begin{pmatrix} 1 & 1 \\ 0 & 1 \\ \end{pmatrix}, \quad \mathrm{and} \quad
h_2 =  \pm \begin{pmatrix} 1 & \omega \\ 0 & 1 \\ \end{pmatrix},
\]
for a fixed element $\omega \in \C - \R$.  Note that the traces of $h_1$ and $h_2$ are determined by the lift to $\SL(2,\C)$, and might not coincide. Thus an arbitrary element $g h_1^m h_2^n \in gH$ can be expressed as
\begin{equation}\label{Eqn:ghForm}
gh_1^m h_2^n \: = \: 
\pm\begin{pmatrix} a & b \\ c & d \\ \end{pmatrix}  \begin{pmatrix} 1 & m + n \omega \\ 0 & 1 \\ \end{pmatrix}
\: = \: \pm \begin{pmatrix} a & a (m + n \omega) + b \\ c &c (m + n \omega) + d  \end{pmatrix} ,
\end{equation}
where the sign $\pm$ depends on the traces of $h_1, h_2$ and the parity of $m,n$.


As above, set  $Z_{\omega} = \{ m + n\omega
\ \vert \ m, n \in \mathbb{Z} \}$ and $Q_{\omega} = \{ m + n \omega \ \vert \ m, n \in \Q \}$. Then
\[
gH \subset  
\left \{ \pm \begin{pmatrix} a & ax + b \\ c & cx + d \\ \end{pmatrix}  \  \Big\vert \  x \in Z_{\omega} \right\}.
\]
Since $g \notin H$, and $H$ is a maximal parabolic subgroup of $\Gamma$, we have $c \neq 0$.
Thus $g h_1^m h_2^n \in gH$  is 
parabolic if and only if $\trace( gh_1^m h_2^n) = \pm(a + d + c(m + n \omega)) \in \{ \pm 2 \}$.  Solving for $x = m + n \omega \in Z_{\omega}$, let
\[
y_+ = y_{+1} = \frac{2 - a - d}{c} \quad \text{and} \quad y_- = y_{-1} = \frac{-2 - a - d}{c}.
\]
  Then the coset $gH$ contains a parabolic element if and only if 
$\{ y_+, y_- \} \cap Z_{\omega} \neq \emptyset$. We will abuse notation slightly by thinking of the subscripts as either symbols ($\pm$) or numbers ($\pm 1$), as convenient.

The pair $y_+$ and $y_-$, corresponding to trace $+2$ and trace $-2$, are the ``problem elements'' that we will need to track throughout the proof. In particular, let $R \subset k$ be the ring generated by the coefficients of the generators of $\Gamma$. 
By expanding $R$ if necessary, we may assume that $c^{-1} \in R$, which implies $y_\pm \in R$.

From here, the proof proceeds as follows. For each number $i \in \{\pm 1\}$, we will construct a homomorphism $\varphi_i \from \Gamma \to G_i$, where $G_i$ is a finite group. Each $\varphi_i$ will be either a congruence quotient of $\Gamma$, or the product of a congruence quotient of $\Gamma$ and a congruence quotient of $\Gamma(\ss_j)$ for some $j$. Then we will package these homomorphisms together to obtain
\[
\varphi = \varphi_- \times \varphi_+ \from \Gamma \longrightarrow G = G_- \times G_+.
\]
In particular, for every $\gamma \in \Gamma$, the image $\varphi(\gamma)$ is a tuple of matrices, each with coefficients in a finite ring, and each with a well-defined trace. To complete the proof, we will see that for every $h \in H$ and every $\ell \in K$, some coordinate of $\varphi(gh)$ differs in trace from the same coordinate of $\varphi(\ell)$. This will imply that $\varphi(gh)$  cannot be conjugate to $\varphi(\ell)$.

For each $i \in \{\pm 1\}$, the definition of $\varphi_i \from \Gamma \to G_i$ depends on whether $y_i$ belongs to $Q_\omega$. 
If $y_i \notin Q_{\omega}$, we argue as follows.

\begin{claim}\label{Claim:EasySep}
Suppose $y_i \in R - Q_{\omega}$. Then there is a ring homomorphism
$\rho_i \from R \rightarrow S_i$, where $S_i$ is a finite ring, such that the following holds. For every pair $(m,n) \in \Z^2$, we have 
\[
\rho_i(m + n \omega) \neq \rho_i(y_i).
\]
\end{claim}

By \Cref{Prop:AddRingSep},
there exist a finite ring $S_i$ and a ring homomorphism
$\rho_i \from R \rightarrow S_i$
such that $\rho_i(y_i) \notin \rho_i(Z_{\omega})$. By the definition of $Z_\omega$, this means $\rho_i(y_i) \neq \rho_i(m + n \omega)$ for any $(m,n)$. \claimqed

Using $\rho_i$, we define a congruence quotient
\begin{equation}\label{Eqn:PhiEasy}
\varphi_i \from \Gamma \hookrightarrow \SL(2, R) \to G_i =  \SL(2, \, S_i) 
\qquad \text{if} \quad  y_i \notin Q_{\omega},
\end{equation}
completing the definition of $G_i$ and $\varphi_i$ in this case.

Alternately, if  $y_i \in Q_{\omega}$,  \Cref{Lem:LinearIndep} provides an infinite set of primes $\PrimeSet$, such that for each $p \in \PrimeSet$ there is an associated finite field $F_\pp$ and ring homomorphism $\eta_p \from R \to F_\pp$. 
The central claim in this case is the following.

\begin{claim}\label{Claim:DehnFillingQuotientSep}
Suppose $y_i \in Q_{\omega}$. Then there is a choice of Dehn filling quotient $\psi_j \from \Gamma \to \Gamma(\ss_j)$ where the coefficients of $\Gamma(\ss_j)$ lie in a finitely generated ring $T_j$, a prime number $p \in \PrimeSet$, and a ring homomorphism $\sigma_{p,i} \from T_j \to E_i$, where $E_i$ is a finite field, such that the following holds. For every pair $(m,n) \in \Z^2$, we have
\[
\eta_p(m + n \omega) \neq \eta_p(y_i) \quad \text{or} \quad
\sigma_{p,i} \circ \trace \circ \,  \psi_j(g h_1^m h_2^n) \neq \pm 2.
\]
\end{claim}

In fact, the desired homomorphism $\sigma_{p,i}$ exists for all sufficiently large $p \in \PrimeSet$. However, we will only need 
$\sigma_{p,i}$ for one $p \in \PrimeSet$. 

We momentarily postpone the proof of \Cref{Claim:DehnFillingQuotientSep} to describe the construction of $\varphi_i$.
The ring homomorphism $\eta_p$ defines a congruence quotient
\[ \nu_0 \from \Gamma \hookrightarrow \SL(2, R) \to \SL(2, F_\pp).
 \]
Similarly, $\sigma_{p,i}$ defines a homomorphism $\nu_i$ factoring through a congruence quotient:
\[
\nu_i \from \Gamma \xrightarrow{ \psi_j  } \psi_j(\Gamma) =  \Gamma(\ss_j) \hookrightarrow \SL(2, T_j) \rightarrow \SL(2, E_i).
\]
We can now define 
\begin{equation}\label{Eqn:PhiHard}
\varphi_i = \nu_0 \times \nu_i \from \Gamma \longrightarrow G_i = \SL(2, F_\pp) \times \SL(2, E_i)
\qquad \text{if} \quad  y_i \in Q_{\omega},
\end{equation}
completing the definition of $G_i$ and $\varphi_i$ in this case.

\medskip

\noindent \emph{Proof of \Cref{Claim:DehnFillingQuotientSep}.}
We begin by specifying the choice of Dehn filling quotient $\psi_1$ or $\psi_2$. Write $y_i = (m_i + n_i \omega)/v_i$ in lowest terms, as in  \Cref{Lem:LinearIndep}. 
If $v_i = 1$, then we set $j = 2$ and work with the Dehn filling $\psi_j = \psi_2 \from \Gamma \to \Gamma(\ss_2)$ for concreteness (although $\psi_1$ would also work).
Assuming $v_i \neq 1$, we have either $v_i \nmid m_i$ or $v_i \nmid n_i$.   If $v_i \nmid m_i$, then we set $j=2$ and select  the Dehn filling $M(\ss_2)$ and the quotient map $\psi_2 \from \Gamma \to \Gamma(\ss_2)$.  Then $\psi_2(H)$
is an infinite cyclic loxodromic subgroup of $\Gamma(\ss_2)$ generated by $\psi_2(h_1)$. Consequently, $\psi_2(h_1^m h_2^n) = \psi_2(h_1^m)$, where
\[
v_i m \not\equiv m_i \smod {{v}_i} \quad \text{because} \quad  v_i \nmid m_i.
\]
Similarly, if $v_i \nmid n_i$, then we set $j=1$ and select the Dehn filling 
 $M(\ss_1)$ and the quotient map $\psi_1 \from \Gamma \to \Gamma(\ss_1)$. This has the effect that $\psi_1(h_1^m h_2^n) = \psi_1(h_2^n)$, where
\[
v_i n \not\equiv n_i \smod {{v}_i} \quad \text{because} \quad v_i \nmid n_i.
\]
In either case, the above non-congruences will be used in the endgame of the proof of the claim.
 Because the arguments for $m_i$ and $n_i$ are entirely parallel, and differ only by a substitution of symbols, we assume without loss of generality that $ v_i \nmid m_i$, hence $j=2$ and we have the Dehn filling quotient $\psi_2 \from \Gamma \to \Gamma(\ss_2)$.

By \Cref{Thm:ThurstonLiftConjugate},
 we can conjugate $\Gamma(\ss_2)$ in $\SL(2, \C)$ 
such that $\Gamma(\ss_2) \subset \SL(2, T)$, where $T = T_2$ is a finitely generated ring in a number field.
Moreover, we may assume that 
\[
\psi_2(g) =  \begin{pmatrix} r & s \\ t & u \\ \end{pmatrix} 
\quad \text{and} \quad 
\psi_2(h_1) =  \begin{pmatrix} \lambda & 0 \\ 0 & \lambda^{-1} \\ \end{pmatrix},
\]
for some $r, s, t, u, \lambda \in \C$ with $\vert \lambda \vert \neq 1$. Then 
\[
\psi_2(gH) = \left\{ \begin{pmatrix} r & s \\ t & u \\ \end{pmatrix}  \begin{pmatrix} \lambda^m & 0 \\ 0 & \lambda^{-m} \\ \end{pmatrix}
\  \Big\vert \  m \in \mathbb{Z}  \right\}  =  \left\{  \begin{pmatrix} r \lambda^m & s  \lambda^{-m}  \\ t \lambda^m & u \lambda^{-m} \end{pmatrix}  \  
\Big\vert \  m \in \mathbb{Z}  \right\},
\]
and $\trace \circ \psi_2( g h_1^m h_2^n) = r \lambda^m + u \lambda^{-m}$.
In particular,  $ \psi_2(gH)$ contains a  parabolic matrix if and only if $r \lambda^m + u \lambda^{-m} = \pm 2$ for some $m$, 
which is true if and only if $\lambda^m$ is a root of $r x^2 \pm 2x + u$.  Let 
\[
\zeta = {\frac{1 + \sqrt{1 - ru}}{r}}, \quad
\xi = {\frac{1 - \sqrt{1 - ru}}{r}}
\]
 be the roots of $rx^2 - 2x + u$.
Then $\{-\zeta, -\xi\}$ are the roots of $rx^2 + 2x + u$. 
By expanding $T = T_2$, if necessary, we may assume that $\{ \zeta, \xi \} \subset T$.
Since $\psi_2(g)$ does not commute with $\psi_2(h_1)$, we have $s \neq 0$ and $t \neq 0$.  Therefore, $ru \neq 1$, which implies that $\zeta \neq \xi$.

To prove the claim, we consider two cases.

\medskip

\noindent
\textbf{Case 1:}  $\lambda^{2m_i} \notin \{ \zeta^{2v_i}, \xi^{2v_i} \}$.

\medskip

By \Cref{Prop:OrderM}, for all sufficiently large primes $p \in \PrimeSet$, 
there exist a finite field $E_i$ and a ring homomorphism $\sigma_{p,i}\from T \rightarrow E_i$, such that 
\[
\sigma_{p,i}(\zeta - \xi) \neq 0, \quad
\sigma_{p,i}(\lambda^{2m_i} - \zeta^{2v_i}) \neq 0, \quad 
\sigma_{p,i}(\lambda^{2m_i} - \xi^{2v_i}) \neq 0,
\]
and the multiplicative order of $\sigma_{p,i}(\lambda)$ is equal to $2p$.  

Suppose for a contradiction that there exist $m, n \in \mathbb{Z}$ such that 
\[ \eta_p(m + n \omega) = \eta_p(y_i) \quad \text{and} \quad
\sigma_{p,i} \circ \trace \circ \psi_2(g h_1^m h_2^n) = \sigma_{p,i}(r \lambda^m + u \lambda^{-m}) = \pm 2.
\]
Then \Cref{Lem:LinearIndep} implies
\[
v_i m \equiv m_i \smod p \quad \Rightarrow \quad 2v_i m \equiv 2 m_i \smod{2p}.
\]
If $\sigma_{p,i}(r \lambda^m + u \lambda^{-m}) = 2$, then $\sigma_{p,i}(\lambda^m)$ is a root of  
$f(x) = \sigma_{p,i}(r) x^2 - 2x + \sigma_{p,i}(u)$ over $E_i$.   Since the two distinct roots of $f$ over $E_i$ are equal to $\sigma_{p,i}(\zeta)$ and $\sigma_{p,i}(\xi)$, we have
$\sigma_{p,i}(\lambda^m) = \sigma_{p,i}(\zeta)$ or $\sigma_{p,i}(\lambda^m) = \sigma_{p,i}(\xi)$.  
Similarly, if $\sigma_{p,i}(r \lambda^m + u \lambda^{-m}) = -2$, then $\sigma_{p,i}(\lambda^m) = \sigma_{p,i}(-\zeta)$ or $\sigma_{p,i}(\lambda^m) = \sigma_{p,i}(-\xi)$.  
In either case, 
\[
\sigma_{p,i}(\lambda^{2v_i m}) = \sigma_{p,i}(\zeta^{2v_i}) \quad \text{or} \quad 
\sigma_{p,i}(\lambda^{2 v_i m}) = \sigma_{p,i}(\xi^{2v_i}).
\]
Since  $2v_i m \equiv 2m_i \smod{2p}$ and the multiplicative order of $\sigma_{p,i}(\lambda)$
is equal to $2p$, this implies that 
\[ 
\sigma_{p,i}(\lambda^{2m_i}) = \sigma_{p,i}(\zeta^{2v_i}) \quad \text{or} \quad 
\sigma_{p,i}(\lambda^{2m_i}) = \sigma_{p,i}(\xi^{2v_i}),
\] 
a contradiction. 

\medskip

\noindent
\textbf{Case 2:} $\lambda^{2m_i} \in \{ \zeta^{2v_i}, \xi^{2v_i} \}$.

\medskip

Without loss of generality, we may assume that $\lambda^{2m_i} = \zeta^{2v_i}$.  
First, suppose that $v_i = 1$. This means that $y_i = m_i + n_i \omega \in Z_\omega$, hence
$g h_1^{m_i} h_2^{n_i} \in gH$ has trace $\pm 2$ and is parabolic in $\Gamma$.  Since $\lambda^{2m_i} = \zeta^{2v_i}$
and $v_i = 1$, $\lambda^{m_i} \in \{\zeta, -\zeta \}$.  Therefore,
$\psi_2(g h_1^{m_i} h_2 ^{n_i}) = \psi_2( g h_1^{m_i})$ is also a parabolic element in $\Gamma(\ss_2)$. 
This is a contradiction, since the only parabolic elements of $\Gamma$ that remain parabolic after the Dehn filling lie in conjugates of $K$,
and the coset $gH$ is disjoint from every conjugate of $K$.  
We conclude that $v_i > 1$.

Next, assume that $\lambda^{2v_i m_i} \neq \xi^{2v_i^2}$. 
By \Cref{Prop:OrderM}, for all sufficiently large primes $p \in \PrimeSet$,
there exist a finite field $E_i$ and a ring homomorphism $\sigma_{p,i}\from T \rightarrow E_i$, such that 
\[
\sigma_{p,i}(\zeta - \xi) \neq 0, \quad
\sigma_{p,i}(\lambda^{2v_i m_i} - \xi^{2v_i^2}) \neq 0,
\]
 and the multiplicative order of $\sigma_{p,i}(\lambda)$ is equal to $2v_ip$.  
Suppose for a contradiction that there exist $m, n \in \mathbb{Z}$ such that 
\[ \eta_p(m + n \omega) = \eta_p(y_i) \quad \text{and} \quad
\sigma_{p,i} \circ \trace \circ \psi_2(g h_1^m h_2^n) = \sigma_{p,i}(r \lambda^m + u \lambda^{-m}) = \pm 2.
\]
Then \Cref{Lem:LinearIndep} implies
\[
v_i m \equiv m_i \smod{p} \quad \Rightarrow \quad
2v_i^2 m \equiv 2v_i m_i \smod{2v_i p}.
\]
Since $\sigma_{p,i}(r \lambda^m + u \lambda^{-m}) = \pm2$, we have
$\sigma_{p,i}(\lambda^m) = \sigma_{p,i}(\pm \zeta)$ or $\sigma_{p,i}(\lambda^m) = \sigma_{p,i}(\pm \xi)$.  
If $\sigma_{p,i}(\lambda^m) = \sigma_{p,i}(\pm\zeta)$, then
$$\sigma_{p,i}(\lambda^{2v_i m}) = \sigma_{p,i}(\zeta^{2v_i}) = \sigma_{p,i}(\lambda^{2m_i}).$$ 
Since the multiplicative order of $\sigma_{p,i}(\lambda)$ is divisible by $2v_i$, we have
\[
2v_i m \equiv 2m_i \smod{2v_i},
\]
contradicting the fact that $ v_i \nmid m_i$.  
If $\sigma_{p,i}(\lambda^m) = \sigma_{p,i}(\pm \xi)$, then, since the multiplicative order of $\sigma_{p,i}(\lambda)$
is equal to $2v_i p$ and $2v_i^2 m \equiv 2 v_i m_i \smod{2v_ip}$, we obtain
$$\sigma_{p,i}(\lambda^{2v_i m_i}) = \sigma_{p,i}(\lambda^{2v_i^2 m}) = \sigma_{p,i}(\xi^{2v_i^2}),$$ 
which contradicts our assumption in choosing $\sigma_{p,i}$.

Finally, assume that $\lambda^{2v_i m_i} = \xi^{2v_i^2}$. 
By \Cref{Prop:OrderM},
there exist a finite field $E_i$ and a ring homomorphism $\sigma_{p,i}\from T \rightarrow E_i$, such that 
$\sigma_{p,i}(\zeta - \xi) \neq 0$, and the multiplicative order of $\sigma_{p,i}(\lambda)$ is divisible by $2v_i^2$.  
Suppose that there exist $m \in \mathbb{Z}$ such that $\sigma_{p,i}(r \lambda^m + u \lambda^{-m}) = \pm 2$.  
Then $\sigma_{p,i}(\lambda^m) = \sigma_{p,i}(\pm \zeta)$ or $\sigma_{p,i}(\lambda^m) = \sigma_{p,i}(\pm \xi)$.  
By the argument above, $\sigma_{p,i}(\lambda^m) \neq \sigma_{p,i}(\pm \zeta)$. 
If $\sigma_{p,i}(\lambda^m) = \sigma_{p,i}(\pm \xi)$, then
$$\sigma_{p,i}(\lambda^{2v_i^2 m}) = \sigma_{p,i}(\xi^{2v_i^2}) = \sigma_{p,i}(\lambda^{2v_i m_i}).$$ Since the multiplicative order of $\sigma_{p,i}(\lambda)$ is divisible by $2v_i^2$,
this implies that 
\[
2v_i^2 m \equiv 2 v_i m_i \smod{2v_i^2}.
\]
  But this contradicts the fact that $ v_i \nmid m_i$, completing the proof of \Cref{Claim:DehnFillingQuotientSep}. 
 \claimqed

We can now complete the proof of the theorem. For each $i \in \{\pm 1\}$, we have defined a homomorphism $\varphi_i \from \Gamma \to G_i$, using \Cref{Eqn:PhiEasy} if $y_i \notin Q_\omega$ and \Cref{Eqn:PhiHard} if $y_i \in Q_\omega$. Now, define
\[
\varphi = \varphi_{-1} \times \varphi_{+1} \from \Gamma \longrightarrow G = G_{-1} \times G_{+1}.
\]
We need to show that $\varphi(K)$ is disjoint from every conjugate of $\varphi(gH)$.

Consider an arbitrary element $gh_1 ^m h_2^n \in gH$, and suppose for a contradiction that $\varphi(gh_1 ^m h_2^n)$ is conjugate 
to $\varphi(\ell)$ for some $\ell \in K$. Since $\ell$ is parabolic, we know $\trace(\ell) \in \{\pm 2\}$. Recalling the general form for 
$gh_1 ^m h_2^n$ in \Cref{Eqn:ghForm}, 
define a number $\epsilon = \epsilon(\ell,m,n) \in \{\pm 1\}$ so that
\[
\trace(g h_1^m h_2^n) = \epsilon \cdot \tfrac{\trace(\ell)}{2} \cdot (a + d + c(m+n \omega)).
\]
In other words, $\epsilon = 1$ when $\trace(h_1^m h_2^n) = \trace(\ell)$, and $\epsilon = -1$ otherwise.
We use the coordinate $\varphi_\epsilon$ of $\varphi$ to obtain a contradiction.

If $y_\epsilon \notin Q_\omega$, then $\varphi_\epsilon(\ell) \in \SL(2,S_\epsilon)$. Since $\varphi_\epsilon(gh_1 ^m h_2^n)$ is conjugate 
to $\varphi_\epsilon(\ell)$, we have
\[
\trace \circ \, \varphi_\epsilon(g h_1^m h_2^n) = \rho_\epsilon \circ \trace(g h_1^m h_2^n) = \rho_\epsilon \left(\epsilon \tfrac{\trace(\ell)}{2} (a + d + c(m+n \omega)) \right) = \rho_\epsilon \big( \trace(\ell) \big)
     = \trace \circ \, \varphi_\epsilon(\ell).
\]
Since $\rho_\epsilon$ is a ring homomorphism and $\epsilon \frac{\trace(\ell)}{2} \in \{ \pm 1 \}$ is a unit, we may rearrange terms to obtain
\[
\rho_\epsilon (a + d + c(m+n \omega))  = \rho_\epsilon \big( \tfrac{2}{\epsilon} \big)
\quad \Rightarrow \quad
\rho_\epsilon (c(m+n \omega))  = \rho_\epsilon \big( \tfrac{2}{\epsilon} - a - d \big).
\]
But then $\rho_\epsilon(m+n \omega) = \rho_\epsilon(y_\epsilon)$, contradicting \Cref{Claim:EasySep}.

If $y_\epsilon \in Q_\omega$, then $\varphi_\epsilon(\ell) = (\nu_0(\ell), \nu_\epsilon(\ell)) \in \SL(2, F_\pp) \times \SL(2,E_{\epsilon})$, a product of two matrices. Since $\nu_0(gh_1 ^m h_2^n)$ is conjugate 
to $\nu_0(\ell)$, we obtain
\[
\trace \circ \, \nu_0(g h_1^m h_2^n) = \eta_p \circ \trace(g h_1^m h_2^n) = \eta_p \left(\epsilon \tfrac{\trace(\ell)}{2} (a + d + c(m+n \omega)) \right) = \eta_p \big( \trace(\ell) \big)
     = \trace \circ \, \nu_0(\ell).
\]
Then the same rearrangement of terms as before gives $\eta_p(m+n \omega) = \eta_p(y_\epsilon)$. Meanwhile, since
$\nu_\epsilon(gh_1 ^m h_2^n)$ is conjugate 
to $\nu_\epsilon(\ell)$, we obtain
\[
\sigma_{p,\epsilon} \circ \trace \circ \,  \psi_j(g h_1^m h_2^n) = \trace \circ \, \nu_\epsilon(g h_1^m h_2^n) =  \trace \circ \, \nu_\epsilon(\ell) = \pm 2,
\]
contradicting \Cref{Claim:DehnFillingQuotientSep}. In either case, the proof is complete.
\end{proof}

\section{Manifolds with non-rectangular cusps}\label{Sec:NonRectangular}

In this section, we prove \Cref{Thm:TrivalentTreeGeneric}. We begin with a cusped, hyperbolic $3$--manifold $M$ containing a horocusp $A$. 
We will construct a sequence of finite covers $\widehat M  \to \mathring M \to M$, with an increasingly strong sequence of properties. The final cover $\widehat M$ will have infinitely many geometric ideal triangulations, implying \Cref{Thm:TrivalentTreeGeneric}. The construction proceeds in four steps. 

\begin{enumerate}[Step 1.]
\item\label{Step:UniqueShortest} Construct a cover $\mathring M \to M$ where $A$ lifts to a horocusp $\mathring A$ that has a unique shortest path to the other cusps. This is accomplished in \Cref{Lem:NearbyCuspsDistinct}.
\smallskip

\item\label{Step:GS} Shrink $\mathring A \subset \mathring M$ to a small sub-horocusp $\mathring{A^t}$. By \Cref{Prop:GS} and \Cref{Cor:EmbeddedAnanas}, the canonical cell decomposition $\mathring \PP$ determined by $\mathring{A^t}$ and the other cusps contains an embedded drilled ananas $\mathring N$ consisting of one or two cells of $\mathring \PP$.
\smallskip

\item\label{Step:LST} Construct a cover $\widehat M \to \mathring M$ where every polyhedron in the lifted polyhedral decomposition $\widehat \PP$ has vertices at distinct cusps. This is accomplished in \Cref{Lem:NoDiagonal}.
\smallskip

\item\label{Step:Subdivide} Now, $\widehat \PP$ can be subdivided into geometric ideal tetrahedra by \Cref{Lem:Preorder}, and furthermore $\widehat \PP$ contains a cover $\widehat N$ of the original ananas $\mathring N$. If the cusp $A$ was non-rectangular, the induced triangulation of $\widehat N$ is equivariant with respect to the cover of $\mathring N$. To conclude the proof, we use \Cref{Lem:Ananas} to find infinitely many ideal triangulations of $\widehat N$, hence of $\widehat M$.
\end{enumerate}

Step \ref{Step:UniqueShortest} builds covers using \Cref{Thm:DoubleCosetAbelian}, whereas
Step \ref{Step:LST} uses \Cref{Prop:PeripheralSeparability}. Steps \ref{Step:GS} and \ref{Step:Subdivide} construct and subdivide polyhedral decompositions, but do not build any covers. The hypothesis that $A$ is a non-rectangular cusp is used only in Step  \ref{Step:Subdivide}. See \Cref{Rem:RectangularProblem} for a detailed description of how this hypothesis is used.

The following basic lemma will be used to apply the results of \Cref{Sec:Separability}.

\begin{lem}\label{Lem:PeripheralCoset}
Let $M = \Hth / \Gamma$ be a cusped hyperbolic manifold. Let $\widetilde B, \widetilde B' \subset \Hth$ be horoballs that cover the same cusp in $M$, and let $g \in \Gamma$ be an isometry such that $g (\widetilde B) = \widetilde B'$. Then the set of \emph{all} elements of $\Gamma$ taking $\widetilde B$ to $\widetilde B'$ is of the form
\[
S =  \Stab_\Gamma(\widetilde B') g \Stab_\Gamma(\widetilde B) = g \Stab_\Gamma(\widetilde B) = \Stab_\Gamma(\widetilde B') \,  g,
\]
both a left coset and a right coset of peripheral subgroups.
\end{lem}

\begin{proof}
Let $h \in S$. Then $h$ differs from $g$ by pre-composition with some element $s \in \Stab_\Gamma(\widetilde B)$ and post-composition with some element $s' \in \Stab_\Gamma(\widetilde B')$. In other words,
\[
S =  \Stab_\Gamma(\widetilde B') g \Stab_\Gamma(\widetilde B) = \big\{ s' \cdot g \cdot s \: \mid \:  s' \in \Stab_\Gamma(\widetilde B'), s \in \Stab_\Gamma(\widetilde B) \big\},
\]
proving the first equality of the lemma. 
Now, observe that $\Stab_\Gamma(\widetilde B') = g \Stab_\Gamma(\widetilde B) g^{-1}$. Thus
\[
S =  g \Stab_\Gamma(\widetilde B) g^{-1} \cdot g \Stab_\Gamma(\widetilde B) = g \Stab_\Gamma(\widetilde B),
\]
proving the second equality of the lemma. The final equality is proved similarly.
\end{proof}

Now, let $M$ be a cusped hyperbolic manifold containing a collection of disjoint, closed horocusps. Let $A$ be one of the horocusps.  Let $\gamma_1, \ldots, \gamma_n$ be the set of all orthogeodesics below a certain length $L$ that connect $A$ to the union of the other cusps. 
This set is finite for any $L$, and non-empty when $L$ is sufficiently large.

Step \ref{Step:UniqueShortest} of the proof is accomplished by the following lemma.

\begin{lem}\label{Lem:NearbyCuspsDistinct}
Let $M = \Hth / \Gamma$ be a cusped hyperbolic manifold and $A \subset M$ a horocusp. For $i = 1, \ldots, n$, let $\gamma_i$ be an  orthogeodesic from $\bdy A$ to some horocusp of $M$. Assume that  $\gamma_i \neq \gamma_j$ for $i \neq j$.
Then there is a finite cover $\mathring f \from \mathring M \to M$,  where $A$ lifts to a horocusp $\mathring A \subset \mathring M$ and where the path-lifts $\mathring \gamma_1, \ldots, \mathring \gamma_n$ that start at $\mathring A$ lead to horocusps of $\mathring M$ that are distinct from one another and from $\mathring A$.
\end{lem}

Recall from \Cref{Def:Lift} that a lift $\mathring A$ must cover $A$ with degree one. Thus, for each $i$, the path $\gamma_i$ has exactly one path-lift $\mathring \gamma_i$ starting at $\mathring A$.

\begin{proof}
Conjugate $\Gamma$ in $\Isom(\Hth)$ so that one preimage of $A$ is a horoball $\widetilde A$ about $\infty$.
The subgroup $K = \Stab_\Gamma(\widetilde A) \cong \Z^2$ can be identified with $\pi_1(A)$. 
Choose a fundamental domain $D \subset \bdy \widetilde A$ for the action of $K$. For each $i$, let $\widetilde \gamma_i$ be a path-lift of $\gamma_i$ to $\Hth$, whose initial point lies in $D$. Let $\widetilde B_i$ be the horoball at the forward endpoint of $\widetilde \gamma_i$. By construction, the orthogeodesics  $\widetilde \gamma_1, \ldots, \widetilde \gamma_n$ lie in distinct $K$--orbits, hence the horoballs $\widetilde B_1, \ldots, \widetilde B_n$ do also. In addition, each $\widetilde B_i$ is disjoint from $\widetilde A$.

For every pair $(i,j)$ with $1 \leq i < j \leq n$, let $S_{ij} \subset \Gamma$ be the set of all deck transformations that map $\widetilde B_j$ to $\widetilde B_i$. This set may be empty (this will be the case if $\widetilde B_i$ and $\widetilde B_j$ cover distinct cusps of $M$). Otherwise, let $g_{ij} \in S_{ij}$ be an arbitrary element and observe that by \Cref{Lem:PeripheralCoset}, $S_{ij} = g_{ij} \Stab_\Gamma(\widetilde B_j)$. Let $T_{ij}$ be the set of all deck transformations that map $\widetilde B_j$ to any horoball in the $K$--orbit of $\widetilde B_i$. If $S_{ij} \neq \emptyset$, we have $T_{ij} = K S_{ij} = K g_{ij} \Stab_\Gamma(\widetilde B_j)$, a double coset of peripheral subgroups of $\Gamma$.

In a similar fashion, let $S_{0j} \subset \Gamma$ be the set of all deck transformations that map $\widetilde B_j$ to $\widetilde A$. If $S_{0j} \neq \emptyset$, \Cref{Lem:PeripheralCoset} says that $S_{0j} = g_{0j} \Stab_\Gamma(\widetilde B_j)$ for an arbitrary element $g_{0j} \in S_{0j}$. Let $T_{0j} = K S_{0j}$.
Finally, define
\[
T = \bigcup_{0 \leq i < j \leq n} T_{ij} = \bigcup_{0 \leq i < j \leq n} K S_{ij} = \bigcup_{i<j, \, S_{ij} \neq \emptyset} K g_{ij} \Stab_\Gamma(\widetilde B_j).
\]

\Cref{Thm:DoubleCosetAbelian} says that for every non-empty $S_{ij}$, the double coset $T_{ij} = K g_{ij} \Stab_{\Gamma}(\widetilde B)$ is separable in $\Gamma$. Thus each such double coset is a closed subset of $ \Gamma$. Since $T$ is a finite union of these closed sets, it follows that $T$ itself is closed, hence separable. Observe that 
$1 \notin T$, because any element $g \in T_{ij}$ moves horoball $\widetilde B_j$ to a distinct location. (This uses the above observation that the horoballs $\widetilde B_i$ and $\widetilde B_j$ lie in distinct $K$--orbits.)

By \Cref{Lem:Separable}, there is a homomorphism $  \varphi \from {  \Gamma} \to   F$, where $  F$ is a finite group, such that $1 =   \varphi(1) \notin \varphi(T)$. Since $  \varphi(K) =   \varphi (K)^{-1}$ is a group, we have
\[
\{ 1 \} \cap   \varphi \left( K \right) \cdot   \varphi \left( \bigcup S_{ij}  \right) = \emptyset
\qquad \Rightarrow \qquad   \varphi \left( K \right) \cap   \varphi \left( \bigcup S_{ij} \right) = \emptyset.
\]
Now, let $ \mathring \Gamma =   \varphi^{-1} \circ   \varphi \big( K \big)$. This is a finite-index subgroup of $  \Gamma$. Let $ \mathring M = \H^3 / { \mathring \Gamma}$. Since $ \Stab_{\Gamma}(\widetilde A) = K =  \Stab_{ \mathring \Gamma}(\widetilde A)$, the horocusp $ \mathring A = \widetilde A / K \subset  \mathring M$ is a lift of $A \subset M$.

By the above displayed equation,
$\varphi( \mathring \Gamma) = \varphi(K)$ is disjoint from $\varphi(S_{ij})$, hence 
$ \mathring \Gamma$ is disjoint from $S_{ij}$ for every $0 \leq i < j \leq n$.  
Thus $\widetilde B_i$ and $\widetilde B_j$ belong to different $ \mathring \Gamma$--orbits and project to distinct cusps in $ \mathring M$. Similarly, $\widetilde A$ and $\widetilde B_j$ belong to different $ \mathring \Gamma$--orbits and project to distinct cusps in $ \mathring M$. Thus the geodesic arcs $ \mathring \gamma_1, \ldots,  \mathring \gamma_n$, namely the quotients of $\widetilde \gamma_1, \ldots, \widetilde \gamma_n$ in $ \mathring M$, lead to cusps of $ \mathring M$ that are distinct from one another and from $ \mathring A$.
\end{proof}

Returning to our plan for proving \Cref{Thm:TrivalentTreeGeneric}, suppose that $A \subset M$ is a horocusp and that  $\{\gamma_1, \ldots, \gamma_n\}$ is the set of all orthogeodesics of minimal length from $A$ to the other cusps (including itself). Following \Cref{Lem:NearbyCuspsDistinct}, we find a cover $\mathring M$ where $A$ lifts to  $\mathring A$ and where the path-lifts $\mathring \gamma_1, \ldots, \mathring \gamma_n$ lead to cusps $\mathring B_1, \ldots, \mathring B_n$ that are distinct from one another and from $\mathring A$. Since the $\mathring B_i$ are distinct, we may adjust their sizes independently. We keep $\mathring B = \mathring B_1$ and shrink $\mathring B_2, \ldots, \mathring B_n$ slightly. Now, $\mathring \gamma_1$ is the unique shortest path in $\mathring M$ from $\mathring A$ to the other horocusps.

In Step \ref{Step:GS}, we apply \Cref{Prop:GS} to shrink $ \mathring A$ by a sufficiently large distance (keeping the name $ \mathring A$) so that  the canonical polyhedral decomposition $ \mathring \PP$ determined by $ \mathring A,  \mathring B$, and the remaining cusps meets $ \mathring A$ in one or two $3$--cells, each with one ideal vertex at $\mathring A$. By \Cref{Cor:EmbeddedAnanas}, these cells glue together to form an embedded drilled ananas $ \mathring N \subset  \mathring M$.

If $A$ is a non-rectangular cusp, then its lift $\mathring A$ is also non-rectangular. Thus, by \Cref{Prop:GS}, the drilled ananas $\mathring N$ consists of two isometric acute-angled tetrahedra.

\begin{define}\label{Def:ReturningDiagonal}
Let $\PP$ be an ideal polyhedral decomposition of a cusped hyperbolic $3$--manifold $M$. A \emph{diagonal} of $\PP$ is a bi-infinite geodesic $\beta$ that is contained in some cell of $\PP$. A diagonal $\beta$ is called \emph{returning} if its endpoints are in the same horocusp of $M$.
\end{define}

Step \ref{Step:LST} of the proof is to apply the following result due to Luo, Schleimer, and Tillmann \cite[Lemmas 8 and 9]{LuoSchTill}.

\begin{lem}\label{Lem:NoDiagonal}
Let $ \mathring M$ be a cusped hyperbolic $3$--manifold with a polyhedral decomposition $ \mathring \PP$. Then there is a finite regular cover $\widehat f \from \widehat M \to  \mathring M$, such that the lifted polyhedral decomposition $\widehat \PP$ has no returning diagonals.
\end{lem}

\begin{proof}
This is proved in \cite[Lemmas 8 and 9]{LuoSchTill}, using a fairly straightforward application of \Cref{Prop:PeripheralSeparability}.
\end{proof}

We can now prove \Cref{Thm:TrivalentTreeGeneric}: a cusped hyperbolic $3$--manifold $M$ containing a non-rectangular cusp $A$ has a cover $\widehat M$ with infinitely many geometric ideal triangulations, organized in a trivalent tree. We stress that each edge of this tree represents a sequence of $n$ geometric Pachner moves as described by \Cref{Cor:AnanasCover}, where $n$ is the degree of the local cover $\widehat{N} \rightarrow \mathring N$ of the ananas $\mathring N\subset \mathring M$.

\begin{proof}[Proof of \Cref{Thm:TrivalentTreeGeneric}]
Let $M$ be a cusped hyperbolic $3$--manifold, and let $A \subset M$ be a non-rectangular cusp. Consider the sequence of finite covers
\[
\widehat M \xrightarrow{\: \widehat f \: } 
\mathring M \xrightarrow{\: \mathring f \: } M
\]
constructed in \Cref{Lem:NearbyCuspsDistinct} and \Cref{Lem:NoDiagonal}. In particular, $ \mathring M$ has a polyhedral decomposition $ \mathring \PP$ such that 
two acute-angled ideal tetrahedra of $ \mathring \PP$ fit together to form a drilled ananas $ \mathring N$. Consequently, $\bdy \mathring N$ consists of two triangular faces of $\mathring \PP$.

Let $\widehat \PP$ be the polyhedral decomposition of $\widehat M$ obtained by pulling back $ \mathring \PP$. Then $\widehat \PP$ contains a submanifold $\widehat N$ that covers $ \mathring N$.
By \Cref{Lem:NoDiagonal}, $\widehat P$ has no returning diagonals, hence the vertices of every polyhedron  $P \subset \widehat \PP$ are mapped to distinct cusps of $\widehat M$.

Choose an ordering $\prec$ on the cusps of $\widehat M$. Then, for every polyhedron  $P \subset \widehat \PP$, we get a total ordering of the vertices of $P$. Thus, by \Cref{Lem:Preorder}, the iterated coning induced by $\prec$ subdivides $\widehat \PP$ into geometric ideal tetrahedra. By construction, $\widehat N$ already consists of tetrahedra, so does not need to be subdivided. 
Now, by \Cref{Cor:AnanasCover}, the initial geometric triangulation of $\widehat N$ (which comes from lifting the two-tetrahedron triangulation of $ \mathring N$) is the start of an infinite sequence of geometric ideal triangulations. 

Since $A$ is non-rectangular, the triangulation of the drilled ananas $ \mathring N$ consists of two acute-angled tetrahedra. Thus every path in the trivalent tree of geometric triangulations of $ \mathring N$ that was described in \Cref{Prop:FareyAnanas} lifts to a path of geometric triangulations of $\widehat N$, hence to a path of geometric triangulations of $\widehat M$.
\end{proof}

\begin{rem}\label{Rem:RectangularProblem}
If $A$ is a rectangular cusp of $M$, the above proof of \Cref{Thm:TrivalentTreeGeneric} still constructs covers $\widehat M \to
\mathring M \to M $, where $\widehat M$ contains a submanifold $\widehat N$ that covers the 
 drilled ananas $\mathring N$. However, this time $\widehat N$ consists of rectangular pyramids that need to be subdivided into tetrahedra. The subdivision imposed by ordering the cusps of $\widehat M$ may impose different choices of diagonals on the rectangles of $\bdy \widehat N$, which would obstruct the triangulation of $\widehat N$ from being equivariant with respect to the cover $\widehat N \to \mathring N$. This means we cannot apply  \Cref{Cor:AnanasCover} to obtain an infinite sequence of triangulations.
The issue of equivariance does not arise if $\widehat N$ is a lift of $\mathring N$.
\end{rem}

In the next section, we will deploy \Cref{Thm:ConjugacySeparability} to construct a cover  $\widehat{M}$ where $\widehat N$ is indeed a lift, enabling us to handle rectangular cusps. This will have the additional benefit that each edge of the trivalent tree will represent a single geometric $2$--$3$ move, as in \Cref{Prop:FareyAnanas}.

\section{Rectangular cusps and Dehn fillings}\label{Sec:RectangularDehn}

In this section, we prove  \Cref{Thm:Main}, which extends \Cref{Thm:TrivalentTreeGeneric} to manifolds with rectangular cusps, and also provides infinitely many geometric triangulations of long Dehn fillings of $\widehat M$.
The proof begins in the same way as Steps \ref{Step:UniqueShortest} and \ref{Step:GS} of
the four-step outline described at the start of \Cref{Sec:NonRectangular}. In particular, we will find a cover $\mathring M \to M$ that  contains a drilled ananas $\mathring N$.
The main challenge, as mentioned in \Cref{Rem:RectangularProblem}, 
is to build further covers where $\mathring N$ continues to lift but (most) returning diagonals stop being returning. Before outlining how to do this, we introduce a definition and a motivating example.

\begin{define}\label{Def:ParabolicEdge}
Let $M = \Hth / \Gamma$ be a cusped hyperbolic $3$--manifold with a distinguished horocusp $A$ and a polyhedral decomposition $\PP$. 
Let $\widetilde A \subset \Hth$ be a horoball covering $A$, and let $\widetilde \PP$ be the lifted polyhedral decomposition. 
An \emph{$\widetilde A$--parabolic diagonal of $\widetilde \PP$} is a bi-infinite geodesic $\widetilde \beta$ contained in a cell of $\widetilde \PP$, whose ends are in horoballs $\widetilde C$ and $\widetilde C'$ such that there is a parabolic isometry $g \in \Stab_\Gamma(\widetilde A)$ with $g (\widetilde C) = \widetilde C'$.

An \emph{$A$--parabolic diagonal of $\PP$} is the projection $\beta \subset M$ of an $\widetilde A$--parabolic diagonal of $\widetilde \PP$, for some horoball $\widetilde A$ covering $A$.
We remark that the choice of $\widetilde A$ is immaterial: for any other horoball $\widetilde A'$ covering $A$, there will be some $\widetilde A'$--parabolic lift of $\beta$. On the other hand, the choice of $\PP$ can affect the collection of $A$--parabolic diagonals, because it affects the set of bi-infinite geodesics that are diagonals of $\PP$ in the first place. We also remark that an $A$--parabolic diagonal is necessarily returning, according to \Cref{Def:ReturningDiagonal}.

A horocusp $C \subset M$ is called \emph{$A$--problematic} (relative to $\PP$) if there is an $A$--parabolic diagonal $\beta$ of $\PP$ whose endpoints lie in $C$. 
\end{define}

\begin{example}\label{Exa:ThornAlwaysProblematic}
Suppose, as in the conclusion of \Cref{Cor:EmbeddedAnanas}, that $M = \Hth/\Gamma$ contains a drilled ananas $N$ that is obtained by gluing one or two cells of $\PP$. Let $A \subset M$ be a horocusp containing the cusp of $N$, and let $B \subset M$ be a horocusp containing the thorn of $N$. 
Now, consider a geodesic $\beta \subset \bdy N$ that lies in a $2$--cell of $\bdy N$. We claim that $\beta$ must be an $A$--parabolic diagonal of $\PP$. Indeed, both endpoints of $\beta$ are in the single thorn of $N$, hence
every lift of $\beta$ to $\Hth$ must have its endpoints in horoballs that  are permuted by $\Stab_\Gamma(\widetilde A)$ for an appropriate horoball $\widetilde A$ covering $A$. Thus, by \Cref{Def:ParabolicEdge}, $B$ is an $A$--problematic cusp relative to $\PP$.

Now, suppose that $\overline M = \Hth/\overline \Gamma$ is some finite cover of $M$ where $A$ lifts to $\overline A$. Since $A$ lifts, or equivalently $\Stab_\Gamma(\widetilde A) = \Stab_{\overline \Gamma}(\widetilde A)$, every $A$--parabolic diagonal $\beta$ of $\PP$ lifts to an $\overline A$--parabolic diagonal $\overline \beta$ of the lifted polyhedral decomposition  $ \overline \PP$.
 In particular, some preimage of $\beta$ continues to be a returning diagonal in $\overline M$.
\end{example}

The gist of the following outline is that \Cref{Exa:ThornAlwaysProblematic} is a worst-case scenario that 
can be isolated and handled.
With \Cref{Thm:ConjugacySeparability} and with enough care, all diagonals that are not in the boundary of a drilled ananas eventually lift to be non-returning, while the ananas continues to lift.

Now, let $M$ be a cusped hyperbolic $3$--manifold containing a horocusp $A$. We will take the following sequence of steps:
\begin{enumerate}[Step 1.]
\item\label{Step:Criterion} Describe a criterion on the distance from $A$ that any $A$--problematic cusp must satisfy. This is accomplished in \Cref{Lem:ProblematicDistance}.

\smallskip
\item\label{Step:AnanasDiagonals} Using the criterion of \Cref{Lem:ProblematicDistance}, find a finite cover $\mathring M \to M$, where $A$ lifts to $\mathring A$. This cover $\mathring M$ contains a polyhedral decomposition $\mathring \PP$ and a drilled ananas $\mathring N$ with its cusp in $\mathring A$ and its thorn in $\mathring B$, such that $\mathring B$ is the only $\mathring A$--problematic cusp of $\mathring M$. In fact, all $\mathring A$--parabolic diagonals of $\mathring \PP$ lie in $\bdy \mathring N$. See \Cref{Lem:RingCoverGeneral} for details.

\smallskip
\item\label{Step:ThornCusp} Using \Cref{Thm:ConjugacySeparability}, find a finite cover $\overline M \to \mathring M$, where all of the above features hold (in particular, $\mathring N$ lifts to $\overline N$), and in addition every polyhedron of $\overline \PP$ has some vertex in a horocusp other than $\overline B$, the thorn of $\overline N$. See \Cref{Lem:NonBlueVertex} for details.

\smallskip
\item\label{Step:FewDiagonals} Using \Cref{Thm:ConjugacySeparability} again, find a finite cover $\widebreve M \to \overline M$, where all of the above features hold (in particular, $\overline N$ lifts to $\widebreve N$), and in addition all returning diagonals have their endpoints in cusps that cover $\overline B$. This means there are very few returning diagonals, and the partial order argument of \Cref{Cor:PreorderTriangulation} suffices to find a geometric triangulation $\widebreve \TT$ that is compatible with infinitely many geometric triangulations of $\widebreve M$. See \Cref{Lem:NonParabolicDiagonal} for details.

\smallskip
\item\label{Step:TwoAnanas} Using $H_1(\widebreve M)$, find a double cover $\widehat M \to \widebreve M$ where $\widebreve N$ has two distinct lifts, called $\widehat N$ and $\widehat N'$. See \Cref{Lem:HomologyCover} for details. We will replace one lift $\widehat N$ with a triangulated solid torus to perform a long Dehn filling $\widehat M(s)$, while using the other lift $\widehat N'$ to obtain infinitely many geometric triangulations of $\widehat M(s)$.
\end{enumerate}

We now proceed to carry out these steps in detail. In the following lemma, $\mathcal{M}(T^2)$ is the moduli space of unit-area flat tori, and $\R_+ \mathcal{M}(T^2)$ is the moduli space of flat tori of any area. 

\begin{lem}\label{Lem:ProblematicDistance}
There is a  function $L \from \R_+ \times \R_+ \mathcal{M}(T^2) \to \R_+$ such that the following holds for every multi-cusped hyperbolic $3$--manifold $M$.

Suppose that $M$ contains a horocusp collection $A, B_1, \ldots, B_k$
and that $\alpha$ is an orthogeodesic
from $A$ to $B_1$.
Then, 
in the canonical polyhedral decomposition $\mathcal \PP$ determined by $A, B_1, \ldots, B_k$, any $A$--problematic horocusp $B_j$ must satisfy $d(A, B_j) < L = L(\len(\alpha), \, \bdy A)$. 

Furthermore, if $A$ is replaced by a sub-horocusp $A^t$ for some $t > 0$, then the distance bound $L$ is replaced by $L+t$. In symbols,
\[
L(\len(\alpha) \!+\! t, \, \bdy A^t) = L(\len(\alpha), \,  \bdy A) + t.
\]
\end{lem}

One particular consequence of \Cref{Lem:ProblematicDistance} is that the length bound $L$ only depends on the horocusps $A$ and $B_1$. Although varying the sizes of $B_2, \ldots, B_k$ may have the effect of changing the polyhedral decomposition $\PP$, thereby changing the collection of $A$--parabolic diagonals of $\PP$, the conclusion of the lemma still holds for the same $L$.

\begin{proof}
Write $M = \Hth/ \Gamma$, and  conjugate $\Gamma$ so that $A$ is covered by a horoball $\widetilde A$ about $\infty$.
Then $\alpha$ has a lift $\widetilde \alpha$ that leads from $\widetilde A$ to a horoball $\widetilde B_1$ covering $B_1$. 
A further conjugation, preserving the point $\infty$, ensures that $\widetilde B_1$ has Euclidean diameter exactly $1$. Then, setting $\ell = \len(\alpha) = \len(\widetilde \alpha)$, it follows that  $\bdy \widetilde A$ lies at Euclidean height $e^\ell$. Let $K = \Stab_\Gamma(\widetilde A)$. 
Note that the length $\ell \in \R_+$ and the Euclidean metric $\bdy A \in \R_+ \mathcal{M}(T^2)$ determine the orbit of horoballs $K\widetilde B_1$ up to Euclidean isometry (equivalently, up to a hyperbolic isometry stabilizing  $\bdy \widetilde A$).

Suppose that $\widetilde \beta \subset \Hth$ is an $\widetilde A$--parabolic diagonal of $\widetilde \PP$.
Let $\widetilde C, \widetilde C'$ be the horoballs at the endpoints of $\widetilde \beta$, and let $\widetilde P$ be a polyhedron containing $\widetilde \beta$.
Then, as in \Cref{Def:CanonicalDecomp}, 
 $\widetilde P$ contains the center of a metric ball $D$ that is tangent to the horoballs about its vertices, including $\widetilde C$ and $\widetilde C'$,
and is disjoint from all other horoballs in the packing.

Let $h$ denote the Euclidean diameter of $\widetilde C$, which 
 is equal to the Euclidean diameter of $\widetilde C'$ because $\widetilde C' \in K \widetilde C$. 
 Then $ d (\widetilde A, \widetilde C) = \ell - \log h $.
 We will see that when $h \ll 1$, or equivalently $d (\widetilde A, \widetilde C) \gg 0$, competing pressures on the diameter of $D$ lead to a contradiction.

Let $w$ be the shortest Euclidean translation length (along $\C = \bdy \Hth \setminus \{\infty\}$) of any element of $K$. Thus the Euclidean distance between the centers of $\widetilde C$ and $\widetilde C'$ is at least $w$. 
 Since $D$ must be tangent to $\widetilde C$ and $\widetilde C'$ but disjoint from $\C$, its Euclidean diameter is bounded below by a function of $w$ and $h$ that grows without bound as $h \to 0$.  On the other hand, a ball of large Euclidean diameter whose lowest point is below Euclidean height $h$ must intersect one of the diameter $1$ horoballs in the $K$--orbit of $\widetilde B_1$. Compare to \Cref{fig:horoball_intro}. Thus every sufficiently small value $h \ll 1$ leads to a contradiction, and there is an upper bound $L$ on $ d (\widetilde A, \widetilde C) = \ell - \log h $. Observe that $L$ depends only on the lattice of horoballs $K \widetilde B_1$, hence on $\ell = \len(\alpha)$ and the Euclidean metric on $\bdy A$. Thus we may write $L = L(\len(\alpha), \bdy A)$.

 To prove the ``furthermore,'' suppose that we replace $A$ by a sub-horocusp $A^t$. This has the effect of replacing $\ell$ by $\ell+t$ and replacing $\widetilde A$ by a horoball $\widetilde A^t$ at Euclidean height $e^{\ell+ t}$. Then the Euclidean length $w$ and the lattice of horoballs $K \widetilde B_1$ both remain the same, hence the same value of $h \ll 1$ leads to a contradiction. However, the distance $d (\widetilde A, \widetilde C)$ has just increased by $t$. Thus replacing $A$ by $A^t$ has the effect of replacing $L$ by $L+t$.
\end{proof}

We can now begin  constructing covers of a cusped hyperbolic $3$--manifold $M$ containing a horocusp $A$. In the following lemma, corresponding to Step \ref{Step:AnanasDiagonals}, we build a cover $\mathring M \to M$ that supports a polyhedral decomposition $\mathring \PP$ that contains
 a drilled ananas $\mathring N$ with its cusp in $\mathring A$, such that all $\mathring A$--parabolic diagonals of $\mathring \PP$ lie in $\bdy \mathring N$. 

\begin{lem}\label{Lem:RingCoverGeneral}
Let $M$ be a cusped hyperbolic $3$--manifold and $A \subset M$ a horocusp. Then there is a finite cover $\mathring f \from \mathring M \to M$ such that the following hold:
\begin{itemize}
\item $\mathring M$ contains a horocusp collection $\mathring A, \mathring B = \mathring B_1, \ldots, \mathring B_k$, where $k \geq 2$ and $\mathring A \subset \mathring M$ is a lift of $A$.
\item There is an orthogeodesic $\mathring \alpha$ from $\mathring A$ to $\mathring B$ that is the unique shortest path from $\mathring A$ to $\cup \mathring B_j$.
\item For large $t > 0$, the polyhedral decomposition $\mathring \PP = \mathring \PP^t$ determined by $\mathring A^t, \mathring B_1, \ldots, \mathring B_k$ contains a drilled ananas $\mathring N$, built out of one or two pyramids whose lateral edges are identified to $\mathring \alpha$, with its cusp in $\mathring A$ and its thorn in $\mathring B$.
\item For every horoball $\widetilde A$ covering $\mathring A$, there is a corresponding preimage $\widetilde N$ of $\mathring N$, such that all $\widetilde A$--parabolic diagonals of $\widetilde \PP$ lie in $\bdy \widetilde N$.
\end{itemize}
\end{lem}

The main point in \Cref{Lem:RingCoverGeneral} is the last bullet, as it provides a partial converse to \Cref{Exa:ThornAlwaysProblematic}.

\begin{proof}
Choose a collection of disjoint horocusps, containing $A$.
Relative to this collection of horocusps, let $\alpha$ be a shortest orthogeodesic in $M$ that starts at $A$. Let $L = L(\len(\alpha), \bdy A)$ be the bound produced by \Cref{Lem:ProblematicDistance}. 
Let $S = \{ \gamma_1, \ldots, \gamma_n \}$ be a set of $n \geq 2$ orthogeodesics starting at $\bdy A$, containing all the orthogeodesics that have length at most $L$. We set $\gamma_1 = \alpha$. By \Cref{Lem:NearbyCuspsDistinct}, there is a finite cover $\mathring M \to M$, where $A$ lifts to a cusp $\mathring A$, and where the $\gamma_i \in S$ have path-lifts $\mathring \gamma_1, \ldots, \mathring \gamma_n$ that start on $\bdy \mathring A$ and lead to cusps that are distinct from one another and from $\mathring A$. Let $\mathring B_i$ be the horocusp at the endpoint of $\mathring \gamma_i$.

Let $\mathring B = \mathring B_1$ be the horocusp at the end of $\mathring \alpha = \mathring \gamma_1$. We keep $\mathring B$ fixed, but shrink each of $\mathring B_2, \ldots \mathring B_n$ to ensure that $d(\mathring A, \mathring B_i) \geq L$ for $i \geq 2, \ldots n$. Any other horocusp of $\mathring M$, labeled $B_i$ for $i = n+1, \ldots, k$, must already satisfy $d(\mathring A, \mathring B_i) \geq L$. In particular, $\mathring \alpha$ is the unique shortest path from $\mathring A$ to 
any cusp of $\mathring M$.

Now, $\mathring M$ and its collection of horocusps satisfies the hypotheses of \Cref{Prop:GS}. Thus, for sufficiently large $t$, we may replace $\mathring A$ with $\mathring A^t$ and build a canonical polyhedral decomposition $\mathring \PP = \mathring \PP^t$ that contains one or two ideal pyramids with a vertex in $\mathring A$ and with their lateral edges glued to $\mathring \alpha$. By \Cref{Cor:EmbeddedAnanas}, these cells glue up to form an embedded, convex drilled ananas $\mathring N \subset \mathring M$. By construction, $\mathring N$ has its cusp in $\mathring A$ and its thorn in $\mathring B$.

By \Cref{Lem:ProblematicDistance}, any $\mathring A$--problematic cusp $\mathring B_i$ in the polyhedral decomposition $\mathring \PP$ must satisfy $d(\mathring A, \mathring B_i) < L$ or equivalently $d(\mathring A^t, \mathring B_i) < L+t$. The only horocusp satisfying these hypotheses is $\mathring B = \mathring B_1$. Thus any $\mathring A$--parabolic diagonal must have its endpoints in $\mathring B$.

Let $\mathring \beta$ be an $\mathring A$--parabolic diagonal of $\mathring \PP = \mathring \PP^t$. In the universal cover $\Hth$, let $\widetilde A$ be a horoball covering $\mathring A$, and let $\widetilde N$ be the preimage of $\mathring N$ containing $\widetilde A^t$. The $\mathring A$--parabolic diagonal $\mathring \beta$ lifts to an $\widetilde A$--parabolic diagonal $\widetilde \beta$, whose endpoints must be in horoballs $\widetilde B, \widetilde B'$ such that $d(\widetilde A, \widetilde B) < L$. By construction, $\mathring \alpha$ is the only orthogeodesic from $\mathring A$ to $\mathring B$ with length less than $L$. Therefore, $\widetilde B$ and $\widetilde B'$ must be
full-sized horoballs connected to $\widetilde A$ by lifts of $\mathring \alpha$. Since the endpoints of $\widetilde \beta$ are in the ideal vertices of $\widetilde N$, and $\widetilde N$ is convex, it follows that $\widetilde \beta \subset \widetilde N$. In particular, $\mathring \beta$ must lie in one of the polyhedra comprising $\mathring N$.

To complete the proof, recall that $\mathring N$ consists of either two ideal tetrahedra or one rectangular-based ideal pyramid. In either case, we think of the constituent cells as pyramids with bases along $\bdy \mathring N$ and lateral edges glued to $\mathring \alpha$. Any diagonal in an ideal pyramid is either a lateral edge or contained in the base. Since the $\mathring A$--parabolic diagonal $\mathring \beta$ has both of its endpoints in $\mathring B$, it cannot be a lateral edge, and must be contained in the base of the ambient pyramid. Thus $\mathring \beta \subset \bdy \mathring N$ and $\widetilde \beta \subset \bdy \widetilde N$, as claimed.
\end{proof}

The next lemma, corresponding to Step \ref{Step:ThornCusp} of the outline, is our first use of  \Cref{Thm:ConjugacySeparability}. Roughly speaking, the lemma says that there is a cover $\overline M \to \mathring M$ where the drilled ananas $\mathring N$ lifts, and where returning diagonals are controlled to a significant degree.

\begin{lem}\label{Lem:NonBlueVertex}
Let $\mathring M = \Hth/ \mathring \Gamma$ be a cusped hyperbolic $3$--manifold containing a distinguished horocusp $\mathring A$ and a drilled ananas $\mathring N$ whose cusp is at $\mathring A$ and whose thorn is in horocusp $\mathring B$. Suppose that $\mathring \PP$ is a polyhedral decomposition of $\mathring M$, with the following property: for every horoball $\widetilde A$ covering $\mathring A$, there is a corresponding preimage $\widetilde N$ of $\mathring N$, such that all $\widetilde A$--parabolic diagonals of $\widetilde \PP$ lie in $\bdy \widetilde N$.
Then
there is a finite cover $\overline f \from \overline M \to \mathring M$ such that the following hold: 
\begin{itemize}
\item $\mathring A$ lifts to a distinguished cusp $\overline A$.
\item $\mathring N$ lifts to a drilled ananas $\overline N$ whose cusp is at $\overline A$ and whose thorn is in horocusp $\overline B$.
\item Every $\overline A$--parabolic diagonal of $\overline \PP$ has its endpoints in $\overline B$.
\item Every polyhedron $\overline P \subset \overline \PP$ has a vertex in some horocusp other than $\overline B$.
\end{itemize}
\end{lem}

\begin{proof}
Let $\widetilde A$ be a horoball covering $\mathring A$, and let $K = \Stab_{\mathring \Gamma}(\widetilde A)$. 
Let $\widetilde N \subset \Hth$ be the preimage of $\mathring N$ containing $\widetilde A$. Then there is a horoball $\widetilde B$ covering $\mathring B$, such that all ideal vertices of $\widetilde N$ lie in $\widetilde A  \cup K \widetilde B$. 
By hypothesis,
all $\widetilde A$--parabolic diagonals of $\widetilde \PP$ lie in $\bdy \widetilde N$ and have their endpoints in horoballs of $K \widetilde B$. 
Let $\mathring P_1, \ldots, \mathring P_n$ be the polyhedra of $\mathring \PP$ that have all of their vertices in $\mathring B$. (If no such polyhedra exist, we may simply set $\overline M = \mathring M$ and let $\overline f$ be the identity map.) For each $\mathring P_i$, let $\mathring \beta_i$ be an edge that is not $\mathring A$--parabolic. Such an edge must exist, because
 all $\mathring A$--parabolic edges belong to $\bdy \mathring N$.

For each $\mathring \beta_i$, choose a preimage $\widetilde \beta_i \subset \Hth$. 
The ends of $\widetilde \beta_i$ lie in horoballs  $\widetilde B_i, \widetilde B'_i$, which must cover $\mathring B$ because both endpoints of $\mathring \beta_i$ are in $\mathring B$.
Thus there is an isometry $g_i \in \mathring \Gamma = \pi_1(\mathring M)$  such that $g_i(\widetilde B_i) = \widetilde B'_i$. By \Cref{Lem:PeripheralCoset}, the set of all isometries in $\mathring \Gamma$ taking $\widetilde B_i$ to $\widetilde B'_i$ is a left coset $g_i \Stab_{\mathring \Gamma}(\widetilde B_i)$. Since $\mathring \beta_i$ is not $\mathring A$--parabolic, the coset $g_i \Stab_{\mathring \Gamma}(\widetilde B_i)$ is disjoint from all $\mathring \Gamma$--conjugates of $K = \Stab_{\mathring \Gamma}(\widetilde A)$. Equivalently, $K$ is disjoint from all $\mathring \Gamma$--conjugates of $g_i \Stab_{\mathring \Gamma}(\widetilde B_i)$.

By \Cref{Thm:ConjugacySeparability}, there is a homomorphism $\varphi_i \from \mathring \Gamma \to G_i$, where $G_i$ is a finite group, such that $\varphi_i(K)$ is disjoint from all $G_i$--conjugates of $\varphi_i \big( g_i \Stab_{\mathring \Gamma}(\widetilde B_i) \big)$. We can now consider the product homomorphism 
\[
\overline \varphi = (\varphi_1, \ldots, \varphi_n) \from \: \mathring \Gamma \longrightarrow G = G_1 \times \dots \times G_n.
\]
Then, for each $i$, the image $\overline \varphi(K)$ is disjoint from all $G$--conjugates of $\overline \varphi \big( g_i \Stab_{\mathring \Gamma}(\widetilde B_i) \big)$.

Now, let $\overline \Gamma = \overline \varphi^{-1} \circ \overline \varphi(K)$, and let $\overline M = \Hth / \overline \Gamma$. We get a covering map $\overline f \from \overline M \to \mathring M$. Then, by construction,  every $\mathring \Gamma$--conjugate of $g_i \Stab_{\mathring \Gamma}(\widetilde B_i)$ is disjoint from $\overline \Gamma$.  %
Since $K = \Stab_{\mathring \Gamma}(\widetilde A) = \Stab_{\overline \Gamma}(\widetilde A)$,
the horocusp $\mathring A \subset \mathring M$ lifts to a horocusp $\overline A = \widetilde A / K \subset \overline M$. Similarly, $\mathring N$ lifts to a drilled ananas $\overline N = \widetilde N / K \subset \overline M$. 
The thorn of $\overline N$ is in the horocusp $\overline B$ that is covered by $\widetilde B$, hence $\overline B$ covers $\mathring B$. This proves the first two bullets in the lemma.

For the next bullet, let $\overline \gamma$ be an $\overline A$--parabolic diagonal in $\overline \PP$. Since $\Hth \to \overline M$ is a regular cover, we may choose a lift $\widetilde \gamma$ that is an $\widetilde A$--parabolic diagonal in $\widetilde \PP$.
Recall that
all $\widetilde A$--parabolic diagonals in $\widetilde \PP$ lie in $\bdy \widetilde N$ and have their endpoints in horoballs of $K \widetilde B$. Since $K \subset \overline \Gamma$, it follows that $\overline \gamma \subset \bdy \overline N$ has its endpoints in $\overline B$, as desired.

To prove the remaining conclusion, let $\overline P$ be a polyhedron of $\overline \PP$. If $\overline P$ has an ideal vertex in some cusp that does not belong to $\overline f^{-1}(\mathring B)$, then certainly $\overline P$ has a vertex that is not in $\overline B$. Otherwise, $\overline P = \overline P_i$ is a lift of some $\mathring P_i$. Let $\overline \beta_i \subset \overline P_i$ be a lift of $\mathring \beta_i \subset \mathring P_i$. We claim that the endpoints of $\overline \beta_i$ are in distinct cusps of $\overline M$, and in particular one endpoint is not in $\overline B$.

Let $\widetilde \gamma_i \subset \Hth$ be an arbitrary preimage of $\overline \beta_i$, and let $\widetilde C_i, \widetilde C'_i$ be horoballs in the packing containing the ends of $\widetilde \gamma_i$. Since $\widetilde \gamma_i$ and $\widetilde \beta_i$ both cover $\mathring \beta_i \subset \mathring M$, there is an isometry $h_i \in \mathring \Gamma$ such that $h_i(\widetilde \gamma_i) = \widetilde \beta_i$, which implies $h_i(\widetilde C_i) = \widetilde B_i$ and $h_i(\widetilde C'_i) = \widetilde B'_i$. Thus the set of all isometries in $\mathring \Gamma$ taking $\widetilde C_i$ to $\widetilde C'_i$ can be written as $h_i^{-1} \cdot  g_i \Stab_{\mathring \Gamma}(\widetilde B_i) \cdot h_i$. By construction,
 this conjugate of $g_i \Stab_{\mathring \Gamma}(\widetilde B_i)$ is disjoint from $\overline \Gamma$. Thus $\widetilde C_i$ and $\widetilde C'_i$ lie in distinct $\overline \Gamma$--orbits, which means that the endpoints of $\overline \beta_i$ 
 are in distinct cusps. This proves the claim and the lemma.
\end{proof}

The next lemma, corresponding to Step \ref{Step:FewDiagonals} of the outline, builds a cover $\widebreve M$ with  even stronger restrictions on the returning diagonals of the polyhedral decomposition $\widebreve \PP$.

\begin{lem}\label{Lem:NonParabolicDiagonal}
Let $\overline M = \Hth/ \overline \Gamma$ be a cusped hyperbolic $3$--manifold containing a distinguished horocusp $\overline A$ and a horocusp $\overline B$. Suppose that $\overline \PP$ is a polyhedral decomposition of $\overline M$, such that every $\overline A$--parabolic diagonal of $\overline \PP$ has its endpoints in $\overline B$.
Then
there is a finite cover $\widebreve f \from \widebreve M \to \overline M$ where $\overline A$ lifts to a distinguished cusp $\widebreve A$, such that every returning diagonal of $\widebreve \PP$ has its endpoints  in $\widebreve f^{-1} (\overline B)$.
\end{lem}

\begin{proof}
Let $R = \{\overline \beta_1, \ldots, \overline \beta_n\}$ be the set of returning diagonals of $\overline \PP$ whose endpoints are not in $\overline B$. This set is finite because $\overline \PP$ has finitely many polyhedra, and each polyhedron contains finitely many diagonals. By hypothesis, every diagonal $\overline \beta_i \in R$ is not $\overline A$--parabolic. Assume that $R \neq \emptyset$, as otherwise we may simply take $\widebreve M = \overline M$.

Let $\widetilde A$ be a horoball covering $\overline A$, and let $K = \Stab_{\overline \Gamma}(\widetilde A)$.
For each $\overline \beta_i$, choose a preimage $\widetilde \beta_i \subset \Hth$. The ends of $\widetilde \beta_i$ lie in horoballs  $\widetilde B_i, \widetilde B'_i$, which cover the same cusp of $\overline M$ because $\overline \beta_i$ is a returning diagonal.
Thus there is an isometry $g_i \in \overline \Gamma = \pi_1(\overline M)$ such that $g_i(\widetilde B_i) = \widetilde B'_i$. By \Cref{Lem:PeripheralCoset}, the set of all isometries in $\overline \Gamma$ taking $\widetilde B_i$ to $\widetilde B'_i$ is a left coset $g_i \Stab_{\overline \Gamma}(\widetilde B_i)$. Since $\overline \beta_i$ is not $\overline A$--parabolic, the coset $g_i \Stab_{\overline \Gamma}(\widetilde B_i)$ is disjoint from all $\overline \Gamma$--conjugates of $K = \Stab_{\overline \Gamma}(\widetilde A)$.

As in the proof of \Cref{Lem:NonBlueVertex}, we use \Cref{Thm:ConjugacySeparability} to find a finite-index subgroup $\widebreve \Gamma \subset \overline \Gamma$ that contains $K$ and is disjoint from every $\overline \Gamma$--conjugate of $g_i \Stab_{\overline \Gamma}(\widetilde B_i)$. Let $\widebreve M = \Hth / \widebreve \Gamma$.
 Since $K \subset \widebreve \Gamma$, the horocusp $\overline A \subset \overline M$ lifts to a horocusp $\widebreve A \subset \widebreve M$. As in the proof of \Cref{Lem:NonBlueVertex}, the disjointness of $\widebreve \Gamma$ and all $\overline \Gamma$--conjugates of $g_i \Stab_{\overline \Gamma}(\widetilde B_i)$ implies that every lift $\widebreve {\beta_i}$ of $\overline \beta_i$ has endpoints in distinct cusps, and is not a returning diagonal. Thus any returning diagonal in $\widebreve M$ must be the preimage of a returning diagonal in $\overline M$ whose endpoints are in $ \overline B$, hence it has endpoints in $\widebreve f^{-1}(\overline B)$.
\end{proof}

Following \Cref{Lem:NonParabolicDiagonal}, the manifold $\widebreve M$ has so few returning diagonals that it is possible to impose a partial order $\prec$ on the cusps of $\widebreve M$ such that every polyhedron $\widebreve P \subset \widebreve \PP$ has a unique lowest vertex. Using \Cref{Cor:PreorderTriangulation}, we can refine $\widebreve \PP$ to a geometric triangulation $\TT$, and apply \Cref{Lem:Ananas} to build infinitely many geometric triangulations of $\widebreve M$. See \Cref{Claim:BreveTriangulation} below for details. 

To find geometric triangulations Dehn fillings, we need to take one more cover.

\begin{lem}\label{Lem:HomologyCover}
Let $\widebreve M = \Hth/ \widebreve \Gamma$ be a cusped hyperbolic $3$--manifold containing at least three cusps and a distinguished horocusp $\widebreve A$. Then there is a double cover $\widehat f \from \widehat M \to \widebreve M$, where $\widebreve A$ has two distinct lifts.
\end{lem}

\begin{proof}
Since $\widebreve M$ is the interior of a compact $3$--manifold with at least three boundary tori, the ``half lives, half dies'' lemma \cite[Lemma 3.5]{Hatcher:3Manifolds}
 implies that $H_1(\widebreve M)$ has a $\Z^n$ direct summand for $n \geq 3$. Let $G$ be the subgroup of $H_1(\widebreve M)$ induced by the inclusion $\widebreve A \to \widebreve M$. Since $G$ has rank at most $2$, there must be a primitive, infinite-order homology class $h \in H_1(\widebreve M)$ such that $\langle h \rangle$ is a direct summand that is linearly independent from $G$. Thus we may define a projection 
$\pi_h \from H_1(\widebreve M) \to \langle h \rangle$ 
such that  $G \subset \ker(\pi_h)$.

Now, consider the sequence of surjective homomorphisms
\[
\widebreve \Gamma = \pi_1(\widebreve M) \xrightarrow{\:  \operatorname{ab} \: } H_1(\widebreve M) \xrightarrow {\:  \pi_h \: } \langle h \rangle \cong \Z \longrightarrow \Z / 2\Z,
\]
where $\operatorname{ab}$ is abelianization. Let $\widehat \varphi \from  \widebreve \Gamma \to \Z/2\Z$ be the composition, and let $\widehat \Gamma = \ker(\widehat \varphi)$. By construction, $G \subset \ker(\pi_h)$, hence $\pi_1(\widebreve A) \subset \ker(\widehat \varphi) = \widehat \Gamma$. By the lifting criterion, $\widebreve A$ lifts to $\widehat M = \Hth / \widehat \Gamma$. Since $\widehat f \from \widehat M \to \widebreve M$ is a regular cover of degree $2$, there must be two distinct lifts.
\end{proof}

We can now prove our main result, \Cref{Thm:Main}.

\begin{proof}[Proof of \Cref{Thm:Main}]
Let $M$ be a cusped hyperbolic $3$--manifold containing a horocusp $A$. Consider the sequence of finite covers
\[
\widehat M \xrightarrow{\: \widehat f \: } 
 \widebreve M \xrightarrow{\: \widebreve f \: } 
\overline M \xrightarrow{\: \overline f \: } 
\mathring M \xrightarrow{\: \mathring f \: } M
\]
constructed in the preceding lemmas. The cusp $A \subset M$ lifts along each covering map. Recall that by \Cref{Lem:RingCoverGeneral}, $ \mathring M$ has at least three cusps, hence $\widebreve M$ does also, and we may indeed apply \Cref{Lem:HomologyCover} to construct $\widehat f$.

By \Cref{Lem:RingCoverGeneral}, $ \mathring M$ has a polyhedral decomposition $\mathring \PP$ such that 
one or two ideal $3$--cells of $ \mathring \PP$ fit together to form a drilled ananas $ \mathring N$ that deformation retracts to $\mathring A$. The $1$--skeleton of $\mathring \PP$ decomposes $\bdy \mathring N$ into two ideal triangles or one ideal rectangle. Furthermore, since $\mathring A$ lifts to $\overline A \subset \overline M$ and $\widebreve A \subset \widebreve M$, the drilled ananas $\mathring N$ lifts to $\overline N \subset \overline M$ and $\widebreve N \subset \widebreve M$.

\begin{claim}\label{Claim:BreveTriangulation}
The polyhedral decomposition $\widebreve \PP$ of $\widebreve M$ can be refined to a geometric ideal triangulation $\widebreve \TT$. If $\bdy \! \widebreve N$ is a single ideal rectangle, then we may choose either diagonal of this rectangle to be an edge in $\widebreve \TT$. Finally, $\widebreve \TT$ is the start of an infinite sequence of geometric triangulations of $\widebreve M$.
\end{claim}

Recall, from \Cref{Lem:NonBlueVertex}, that the drilled ananas $\overline N \subset \overline M$ has its thorn in a horocusp $\overline B$, such that every polyhedron $\overline P \subset \overline \PP$ has a vertex in some horocusp apart from $\overline B$. Thus, in the lifted polyhedral decomposition $\widebreve \PP$ of $\widebreve M$, every polyhedron $\widebreve P$ must have at least one vertex in a horocusp that is not in $\widebreve f^{-1}(\overline B)$. We call the cusps of $\widebreve f^{-1}(\overline B)$ \emph{blue}. By \Cref{Lem:NonParabolicDiagonal}, all returning diagonals of $\widebreve \PP$ have their endpoints in blue cusps.

Let $V$ be the set of cusps of $\widebreve M$. We impose a partial order $\prec$ on $V$ as follows:  the non-blue cusps are totally ordered in some fashion; the blue cusps are pairwise incomparable; and $\widebreve C \prec \widebreve B$ for every blue cusp $\widebreve B$ and non-blue cusp $\widebreve C$. Since every polyhedron  $P \subset \widebreve \PP$ has at least one non-blue ideal vertex, and the non-blue vertices of $P$ are totally ordered below the blue ones, it follows that $P$ has a unique $\prec$--minimal vertex. Thus, by \Cref{Lem:Preorder}, the iterated coning of $\PP$ induced by $\prec$ produces $\widebreve {\PP'}$, a well-defined subdivision of $\widebreve \PP$ into geometric ideal pyramids. By \Cref{Cor:PreorderTriangulation}, any choice of diagonals in the non-triangular faces of $\widebreve {\PP'}$ produces a geometric ideal triangulation $\widebreve \TT$. 

By construction, the thorn of $\widebreve N$ is in a blue cusp of $\widebreve f^{-1}(\overline B)$. Thus the partial order $\prec$ does not impose any ordering on the ideal vertices of $\bdy \! \widebreve N$. If this boundary is a single ideal rectangle, the coning induced by $\prec$ in \Cref{Lem:Preorder} does not subdivide it, and we may choose our preferred diagonal to subdivide $\widebreve N$ into two ideal tetrahedra. By \Cref{Lem:Ananas},  $\widebreve N$ admits infinitely many geometric ideal triangulations, hence $\widebreve M$ does also. This proves the claim. \claimqed

\smallskip

Next, we lift the triangulation $\widebreve \TT$ of $\widebreve M$ to a geometric triangulation $\widehat \TT$ of $\widehat M$. Then the infinite sequence of geometric triangulations of $\widebreve M$ lifts to an infinite sequence of geometric triangulations of $\widehat M$, as claimed in the statement of the theorem. 

By \Cref{Lem:HomologyCover}, the horocusp $\widebreve A \subset \widebreve M$ has two distinct lifts to $\widehat M$, which we call $\widehat A$ and $\widehat A'$. Consequently, the drilled ananas $\widebreve N \supset \widebreve A$ also has two distinct lifts to $\widehat M$, namely $\widehat N \supset \widehat A$ and $\widehat N' \supset \widehat A'$. The two distinct lifts of $\widebreve A$ and $\widebreve N$ play distinct roles in the Dehn filling argument.

\begin{claim}\label{Claim:SpunTriangulation}
For all but finitely many choices of slope $s$ on $\bdy \widehat A$, the following hold:
\begin{itemize}
\item The Dehn filled manifold $\widehat M(s)$ has a hyperbolic structure where the core curve  $\gamma$ of the filled solid torus is isotopic to a closed geodesic $\gamma_s$.
\item The ideal tetrahedra of $\widehat \TT$ remain geometric in the hyperbolic metric on $\widehat M(s)$. Each ideal vertex of a tetrahedron of $\widehat \TT$ that used to enter cusp  $\widehat A$ now spins about the geodesic $\gamma_s$.
\item With an appropriate choice of diagonals in \Cref{Claim:BreveTriangulation}, the union of the two tetrahedra in $\widehat N$ has convex boundary in $\widehat M(s)$.
\end{itemize}
\end{claim}

The first two bullets in the claim follow from Thurston's hyperbolic Dehn surgery theorem \cite[Chapter 4]{Thurston:Notes}. 
For a sufficiently long slope $s$, the hyperbolic metric on $\widehat M(s)$ is obtained via an arbitrarily small deformation of the metric on $\widehat M$. Thus, for every ideal tetrahedron $T \subset \widehat \TT$, a sufficiently small deformation of the metric will keep $T$ geometric and positively oriented. As in \cite[Section 4.4]{Thurston:Notes}, the two tips of ideal tetrahedra that enter $\widehat A$ will now spin about the core geodesic $\gamma_s$.

For the last bullet of the claim, suppose first that the cusp $A \subset M$ is non-rectangular (hence, so are its lifts). Then \Cref{Cor:EmbeddedAnanas} implies that the original drilled ananas $\mathring N$ consists of acute tetrahedra, and $\bdy \mathring N$ is strictly convex at all three of its edges. The same properties are preserved in the lift $\widehat N$ and are still preserved in $\widehat N(s)$ after a sufficiently small deformation of the metric. 

Next, suppose that $A \subset M$ is rectangular. Then \Cref{Cor:EmbeddedAnanas} implies that the original drilled ananas $\mathring N$ has convex boundary, with interior angles strictly less than $\pi$ at the two edges of $\bdy \mathring N \cap \mathring \PP$, and an angle of $\pi$ along the (arbitrary) diagonal of the ideal rectangle. The same properties remain true in the lift $\widehat N \subset \widehat M$. When we deform the metric on $\widehat M$ to obtain $\widehat M(s)$, the interior angle of $\pi$ may become $\pi + \epsilon$ for small $\epsilon$, violating convexity, but then the opposite choice of diagonal on $\bdy \widehat N$ will have interior angle $\pi - \epsilon$. Thus an appropriate choice of diagonal in \Cref{Claim:BreveTriangulation} keeps $\bdy \widehat N$ convex in $\widehat M(s)$. \claimqed

To construct geometric triangulations of $\widehat M(s)$, we need to introduce the solid-torus analogue of a drilled ananas. 
A \emph{filled ananas} is a $3$--manifold $X$ homeomorphic to a solid torus with one boundary point removed, and endowed with a complete hyperbolic metric with the following properties. The boundary $\bound X$ is made up of two totally geodesic ideal triangles, with vertices at the removed point. These two ideal triangles are glued by isometry along their edges to form a standard two-triangle triangulation of a once-punctured torus, with shearing and bending allowed along the edges. Furthermore, $X$ is subdivided into geometric ideal tetrahedra.

\begin{claim}\label{Claim:FilledAnanas}
For all but finitely many choices of slope $s$ on $\bdy \widehat A$, the hyperbolic manifold $\widehat M(s)$ has a geometric triangulation $\widehat \TT(s)$ with the following properties. Finitely many tetrahedra of $\widehat \TT(s)$ fit together to form a filled ananas $X(s)$. 
Furthermore, the restriction of $\widehat \TT(s)$ to the complement $\widehat M(s) \setminus X(s)$  is combinatorially isomorphic to the restriction of $\widehat \TT$ to the complement $\widehat M \setminus \widehat N$.
\end{claim}

This statement is due to Gu\'eritaud and Schleimer, and closely resembles \cite[Theorem 1]{GS:canonical}. Assuming that $\widehat \TT$ is the canonical triangulation of $\widehat M$ with respect to some choice of horocusps, they construct the triangulated filled ananas $X(s)$ and endow it with a geometric structure isometric to the completion of the two spun tetrahedra of $\widehat N$ mentioned in \Cref{Claim:SpunTriangulation}. Then they replace $\widehat N$ with $X(s)$, and prove
 that the resulting triangulation $\widehat \TT(s)$ is the canonical triangulation of $\widehat M(s)$. In fact, the construction of $X(s)$, which occurs in \cite[Section 2]{GS:canonical} and is encapsulated in \cite[Corollary 16]{GS:canonical}, only uses the hypotheses that $\widehat \TT$ is geometric and that $\widehat N$ remains convex in $\widehat M(s)$. See also \cite[Corollary 4.18]{HamPurcell2020}. The combinatorial structure of the triangulation of $X(s)$ is closely guided by the combinatorics of the Farey graph $\mathcal{F}$ and the continued fraction expansion of the filling slope $s$, while the hyperbolic metric on $X(s)$ is constructed using Casson and Rivin's work on angle structures and volume optimization \cite{Rivin:Volume, FG:AngledSurvey}. In particular, the canonicity of $\widehat \TT$ is not needed in the proof that $\widehat \TT(s)$ is geometric. \claimqed

To complete the proof of the theorem, we have

\begin{claim}
For all but finitely many choices of slope $s$ on $\bdy \widehat A$, the Dehn filled manifold $\widehat M(s)$ has an infinite sequence of geometric triangulations connected by geometric $2$--$3$ moves.
\end{claim}

Observe that by \Cref{Claim:FilledAnanas}, the geometric triangulation $\widehat \TT(s)$ agrees with $\widehat \TT$ on the drilled ananas $\widehat N' \subset \widehat M(s)$. By \Cref{Lem:Ananas}, this two-tetrahedron geometric triangulation of $\widehat N'$ is the start of an infinite sequence of geometric triangulations connected by geometric $2$--$3$ moves.
\qquad \qquad \FourStarOpen 
\end{proof}

The following remark states a version of \Cref{Thm:TrivalentTreeGeneric} for manifolds with rectangular cusps.
In the statement, a \emph{$4$--$4$ move}
is a local move on triangulations, which takes an octahedron that has been decomposed into 4 tetrahedra along one of its three internal diagonals and replaces it with a decomposition into 4 tetrahedra along a different internal diagonal.
The move is called a \emph{geometric $4$--$4$ move} if both decompositions are into geometric ideal tetrahedra.
\begin{rem}\label{Rem:RectangularFarey}
If $A \subset M$ is a rectangular cusp, the cover $\widebreve M$ contains an infinite trivalent tree of geometric ideal triangulations, where one edge of the tree is a geometric $4$--$4$ move and the remaining edges are geometric $2$--$3$ moves.
This can be seen as follows. 
In \Cref{Claim:BreveTriangulation} of the above proof, the drilled ananas $\widebreve N$ consists of an ideal rectangular pyramid $\widebreve P \subset \widebreve \PP$. After the pyramidal decomposition induced by $\prec$, the $3$--cell $\widebreve P'$ glued to $\widebreve P$ is also an ideal rectangular pyramid. The two choices of diagonal for the shared face of $\widebreve P \cap \widebreve P'$ lead to ideal triangulations that differ by a geometric $4$--$4$ move. Each of these choices can serve  in \Cref{Lem:Ananas} as the starting configuration in an infinite sequence of geometric ideal triangulations. In the proof of \Cref{Prop:FareyAnanas}, the dual tree of the Farey complex splits in half along the edge $\tfrac{0}{1}-\tfrac{1}{0}$: half of the tree is reachable by geometric $2$--$3$ moves if we choose the diagonal $\tfrac{1}{1}$, and the other half is reachable if we choose the diagonal $\frac{-1}{1}$. The union of these halves is an infinite trivalent tree of geometric $2$--$3$ moves, with one edge replaced by a  $4$--$4$ move.
\end{rem}

We close the paper by pointing out that the figure--8 knot complement $M$ does not contain a drilled ananas, because $M$ has only one cusp. Nevertheless, by the work of Dadd and Duan \cite{DaddDuan2016}, this manifold has an infinite family of geometric triangulations. 
Thus the presence of a drilled ananas is a sufficient but not necessary condition. We wonder what other features will also guarantee an infinite sequence of geometric triangulations.

\bibliographystyle{plain}
\bibliography{infinite_bib}

\begin{thebibliography}{10}

\bibitem{akiyoshi_finiteness}
Hirotaka Akiyoshi.
\newblock Finiteness of polyhedral decompositions of cusped hyperbolic
  manifolds obtained by the {E}pstein-{P}enner's method.
\newblock {\em Proc. Amer. Math. Soc.}, 129(8):2431--2439, 2001.

\bibitem{Baker-Cooper:QFS}
Mark~D. Baker and Daryl Cooper.
\newblock Finite-volume hyperbolic 3-manifolds contain immersed
  quasi-{F}uchsian surfaces.
\newblock {\em Algebr. Geom. Topol.}, 15(2):1199--1228, 2015.

\bibitem{bass1980groups}
Hyman Bass.
\newblock Groups of integral representation type.
\newblock {\em Pacific J. Math.}, 86(1):15--51, 1980.

\bibitem{ChagasZalesskii}
Sheila~C. Chagas and Pavel~A. Zalesskii.
\newblock Hyperbolic 3-manifolds groups are subgroup conjugacy separable.
\newblock {\em Ann. Sc. Norm. Super. Pisa Cl. Sci. (5)}, 19(3):951--963, 2019.

\bibitem{ChampanerkarKofmanPurcell}
Abhijit Champanerkar, Ilya Kofman, and Jessica~S. Purcell.
\newblock Geometry of biperiodic alternating links.
\newblock {\em J. Lond. Math. Soc. (2)}, 99(3):807--830, 2019.

\bibitem{Culler:Lifting}
Marc Culler.
\newblock Lifting representations to covering groups.
\newblock {\em Adv. in Math.}, 59(1):64--70, 1986.

\bibitem{SnapPy}
Marc Culler, Nathan~M. Dunfield, Matthias Goerner, and Jeffrey~R. Weeks.
\newblock Snap{P}y, a computer program for studying the geometry and topology
  of $3$-manifolds.
\newblock Available at \url{http://snappy.computop.org}.

\bibitem{CuSh}
Marc Culler and Peter~B. Shalen.
\newblock Varieties of group representations and splittings of {$3$}-manifolds.
\newblock {\em Ann. of Math. (2)}, 117(1):109--146, 1983.

\bibitem{DaddDuan2016}
Blake Dadd and Aochen Duan.
\newblock Constructing infinitely many geometric triangulations of the figure
  eight knot complement.
\newblock {\em Proc. Amer. Math. Soc.}, 144(10):4545--4555, 2016.

\bibitem{EpsteinPenner}
David B.~A. Epstein and Robert~C. Penner.
\newblock Euclidean decompositions of noncompact hyperbolic manifolds.
\newblock {\em J. Differential Geom.}, 27(1):67--80, 1988.

\bibitem{FG:AngledSurvey}
David Futer and Fran\c{c}ois Gu\'eritaud.
\newblock From angled triangulations to hyperbolic structures.
\newblock In {\em Interactions between hyperbolic geometry, quantum topology
  and number theory}, volume 541 of {\em Contemp. Math.}, pages 159--182. Amer.
  Math. Soc., Providence, RI, 2011.

\bibitem{goerner2017geodesic}
Matthias Goerner.
\newblock Geodesic triangulations exist for cusped {P}latonic manifolds.
\newblock {\em New York J. Math.}, 23:1363--1367, 2017.

\bibitem{GueriFut}
Fran\c{c}ois Gu\'eritaud.
\newblock On canonical triangulations of once-punctured torus bundles and
  two-bridge link complements.
\newblock {\em Geom. Topol.}, 10:1239--1284, 2006.
\newblock With an appendix by David Futer.

\bibitem{Gueritaud:SolidTori}
Fran\c{c}ois Gu\'eritaud.
\newblock Deforming ideal solid tori.
\newblock arXiv:0911.3067, 2009.

\bibitem{GS:canonical}
Fran\c{c}ois Gu\'{e}ritaud and Saul Schleimer.
\newblock Canonical triangulations of {D}ehn fillings.
\newblock {\em Geom. Topol.}, 14(1):193--242, 2010.

\bibitem{HamPurcell2020}
Sophie~L. Ham and Jessica~S. Purcell.
\newblock Geometric triangulations and highly twisted links.
\newblock arXiv:2005.11899, 2020.

\bibitem{hamilton2005finite}
Emily Hamilton.
\newblock Finite quotients of rings and applications to subgroup separability
  of linear groups.
\newblock {\em Trans. Amer. Math. Soc.}, 357(5):1995--2006, 2005.

\bibitem{HWZ:Separability}
Emily Hamilton, Henry Wilton, and Pavel~A. Zalesskii.
\newblock Separability of double cosets and conjugacy classes in 3-manifold
  groups.
\newblock {\em J. Lond. Math. Soc. (2)}, 87(1):269--288, 2013.

\bibitem{Hatcher:3Manifolds}
Allen Hatcher.
\newblock Notes on basic 3-manifold topology.
\newblock \url{https://pi.math.cornell.edu/~hatcher/3M/3Mdownloads.html}, 2007.

\bibitem{HIKMOT16}
Neil Hoffman, Kazuhiro Ichihara, Masahide Kashiwagi, Hidetoshi Masai, Shin'ichi
  Oishi, and Akitoshi Takayasu.
\newblock Verified computations for hyperbolic 3-manifolds.
\newblock {\em Experimental Mathematics}, 25(1):66--78, 2016.

\bibitem{Long:TotallyGeodesic}
Darren~D. Long.
\newblock Immersions and embeddings of totally geodesic surfaces.
\newblock {\em Bull. London Math. Soc.}, 19(5):481--484, 1987.

\bibitem{LuoSchTill}
Feng Luo, Saul Schleimer, and Stephan Tillmann.
\newblock Geodesic ideal triangulations exist virtually.
\newblock {\em Proc. Amer. Math. Soc.}, 136(7):2625--2630, 2008.

\bibitem{NeumannZagier}
Walter~D. Neumann and Don Zagier.
\newblock Volumes of hyperbolic three-manifolds.
\newblock {\em Topology}, 24(3):307--332, 1985.

\bibitem{Nimershiem}
Barbara~E. Nimershiem.
\newblock Isometry classes of flat {$2$}-tori appearing as cusps of hyperbolic
  {$3$}-manifolds are dense in the moduli space of the torus.
\newblock In {\em Low-dimensional topology ({K}noxville, {TN}, 1992)}, Conf.
  Proc. Lecture Notes Geom. Topology, III, pages 133--142. Int. Press,
  Cambridge, MA, 1994.

\bibitem{Petronio:IdealTriangulations}
Carlo Petronio.
\newblock Ideal triangulations of hyperbolic {$3$}-manifolds.
\newblock {\em Boll. Unione Mat. Ital. Sez. B Artic. Ric. Mat. (8)},
  3(3):657--672, 2000.

\bibitem{PetronioPorti}
Carlo Petronio and Joan Porti.
\newblock Negatively oriented ideal triangulations and a proof of {T}hurston's
  hyperbolic {D}ehn filling theorem.
\newblock {\em Expo. Math.}, 18(1):1--35, 2000.

\bibitem{Rivin:Volume}
Igor Rivin.
\newblock Euclidean structures on simplicial surfaces and hyperbolic volume.
\newblock {\em Ann. of Math. (2)}, 139(3):553--580, 1994.

\bibitem{Sirotkina}
M.~L. Sirotkina.
\newblock On the triangulation of three-dimensional hyperbolic manifolds.
\newblock {\em Algebra i Analiz}, 14(6):192--204, 2002.

\bibitem{Thurston:Notes}
William~P. Thurston.
\newblock The geometry and topology of three-manifolds.
\newblock \url{http://library.msri.org/books/gt3m/}, 1980.

\bibitem{wada1996inequality}
Masaaki Wada, Yasushi Yamashita, and Han Yoshida.
\newblock An inequality for polyhedra and ideal triangulations of cusped
  hyperbolic {$3$}-manifolds.
\newblock {\em Proc. Amer. Math. Soc.}, 124(12):3905--3911, 1996.

\bibitem{WZ:DistinguishingGeometries}
Henry Wilton and Pavel Zalesskii.
\newblock Distinguishing geometries using finite quotients.
\newblock {\em Geom. Topol.}, 21(1):345--384, 2017.

\bibitem{yoshida1996ideal}
Han Yoshida.
\newblock Ideal tetrahedral decompositions of hyperbolic {$3$}-manifolds.
\newblock {\em Osaka J. Math.}, 33(1):37--46, 1996.

\end{thebibliography}

\end{document}